\pdfoutput=1
\documentclass[11pt]{article}
\usepackage[utf8]{inputenc}
\usepackage[a4paper, margin=1in]{geometry}
\usepackage{amsmath, amssymb, amsthm, mathtools}
\usepackage{graphicx, caption, subcaption}
\usepackage{enumitem}
\usepackage{hyperref}
\usepackage{color}
\usepackage{float}
\usepackage{xcolor}
\usepackage{natbib}

\usepackage{times}

\usepackage{setspace}
\usepackage{changepage}

\title{\textbf{Some stability results for a Fractional Differential Equation with two delays}}
\author{  Pragati Dutta, Sachin Bhalekar*\\
	\textit{School of Mathematics and Statistics, University of Hyderabad,India}\\
	*corresponding author email id:sachinbhalekar@uohyd.ac.in}
\date{}

\newtheorem{theorem}{Theorem}[section]

\theoremstyle{remark}

\begin{document}
	\maketitle

\begin{abstract}
We investigate a nonlinear scalar Caputo fractional delay differential equation with two discrete delays and a delay-dependent feedback coefficient. Under standard Lipschitz assumptions, existence and uniqueness of solutions are established through a fixed-point argument and a fractional Gronwall inequality. The trivial equilibrium is then studied by linearization and characteristic-root analysis. When the first delay is zero, explicit delay-independent stability and instability regions are obtained, together with critical-delay conditions that account for the varying coefficient. When both delays are retained, sufficient parameter conditions for delay-independent stability and instability are derived, and an explicit sufficient threshold for the existence of a positive real characteristic root is obtained. Possible Hopf boundaries are identified through purely imaginary roots and the associated transversality condition. The influence of the fractional order on the spectral conditions is highlighted, and the analytical results are illustrated by numerical simulations and stability diagrams.
\end{abstract}


\section{Introduction}

In 1695, Leibnitz proposed a notation $\frac{d^nx}{dt^n}$ for the $n^{th}$ order derivative of a function $x(t)$, where $n$ is a natural number. In his letter dated \date{30 September 1695}, L’Hopital \cite{leibniz1860leibnizens,podlubny1998fractional} asked the meaning of this term when $n=1/2$. This is popularly known as the birth of the fractional calculus, a branch of mathematics that studies operators of noninteger order. By Fractional Derivative (FD), one means the derivative of an arbitrary order, e.g, a real or complex number, a function of t, etc.. The fractional derivative operators are non-local, unlike the classic integer-order derivative \cite{bhalekar2019can,li2007remarks}. Thus, one must provide the state values from the initial point when evaluating the FDs. This procedure proved extremely useful while modeling the memory properties in natural systems \cite{Podlubny1999,kilbas2006theory}. This also gives the flexibility of selecting an order that fits the experimental data \cite{mainardi2022fractional,bagley1983theoretical}.\\
Similarly, the Delay Differential Equations (DDE) contain the states at the past time values \cite{smith2010introduction,lakshmanan2011dynamics}. Though the DDEs look very similar to their non-delayed counterparts, they are extremely difficult to analyze analytically \cite{hale2013introduction}. A single scalar DDE can have oscillatory solutions. Moreover, one can have a scalar DDE exhibiting chaos \cite{sprott2007simple,senthilkumar2005bifurcations,ruiz2013chaos}. The DDEs are infinite-dimensional dynamical systems that occur in various models in physics \cite{otto2019nonlinear,bazighifan2021oscillation}, chemistry \cite{epstein1991differential,roussel1996use}, biology and medicine \cite{rihan2021delay,bani2017analysis,glass2021nonlinear,dong2023application}, engineering and control systems \cite{insperger2011semi,michiels2010control,gumussoy2014computer}, finance \cite{nagarajan2025stock,flore2019feynman,agrawal2020jump}, and so on.\\
The fractional order differential equations with delay(FDDE) contain both these tools: the fractional derivative as well as the delay. Thus, one gets a better model for the given natural system. The numerical methods to solve FDDEs are presented in \cite{bhalekar2011predictor,daftardar2014new}. The stability and bifurcation analysis of these systems is developed in \cite{bhalekar2016stability,bhalekar2011fractional,gupta2024fractional,bhalekar2022stability}. The examples of chaotic FDDEs are \cite{bhalekar2012dynamical,daftardar2012dynamics,bhalekar2012generalized}. The DDEs involving multiple delays occur in applications involving two or more processes \cite{gu2003stability,belair1995age,niculescu2002delay}. The logistic model with two delays is studied in \cite{braddock1983two}. Piotrowska proposed a model with two delays in tumor growth \cite{piotrowska2008hopf}. Boullu et al employed such equations in the production of platelets \cite{boullu2020stability}.\\
There are very few papers devoted to the analysis of such equations. Hale and Huang \cite{hale1993global}, Belair and Campbell \cite{belair1994stability}, Li, Ruan, and Wei \cite{li1999stability}, and Bhalekar \cite{bhalekar2019analysing} studied particular classes of DDEs with multiple delays. Bhalekar and Dutta \cite{bhalekar2025analysis} discussed the stability of a class of FDDEs with two delays.

{
In this paper, we study the stability results for nonlinear Caputo fractional delay differential equation
\[
D^\alpha x(t)=-\gamma x(t)+g\bigl(x(t-\tau_1)\bigr)
-e^{-\gamma\tau_2}g\bigl(x(t-\tau_1-\tau_2)\bigr),
\qquad 0<\alpha\leq 1
\]
near the trivial equilibrium point.

The principal mathematical difficulty in the present problem is that the characteristic equation contains both the fractional power $\lambda^\alpha$ and two coupled exponential terms involving $\tau_1$ and $\tau_1+\tau_2$. Moreover, the gain of the second delayed term, $e^{-\gamma\tau_2}$, depends on the delay itself. Consequently, changing $\tau_2$ alters both the phase and the magnitude of the feedback. Standard single-delay stability criteria, which assume constant coefficients, therefore cannot be applied directly, and the possible crossings of characteristic roots through the imaginary axis cannot generally be reduced to a polynomial condition.

We overcome this difficulty in two stages. When $\tau_1=0$, the equation is recast as a single-delay FDDE with the moving coefficient $b(\tau)=-k e^{-\gamma\tau}$; tracking this coefficient through the known stability regions yields explicit delay-independent regions and critical-delay conditions. When both delays are retained, we combine the real and imaginary parts of the characteristic equation with uniform trigonometric bounds to exclude imaginary-axis crossings in specified parameter regions. Positive-real-root and extremum arguments are then used to establish delay-independent instability and an explicit feedback threshold $k_*$. Thus, the contribution is a set of analytically verifiable stability and instability regions that remain valid uniformly over the relevant delays and fractional orders, despite the delay-dependent coefficient.

More specifically, Deng, Li, and Lü \cite{deng2007stability} established a general characteristic-root criterion for linear fractional systems with multiple constant delays, but did not derive explicit parameter regions for a nonlinear scalar equation whose feedback coefficient varies with a delay. The two-delay fractional equation studied by Bhalekar and Dutta \cite{bhalekar2025analysis} has the form $D^\alpha x(t)=a x(t)+b x(t-\tau)-b x(t-2\tau)$; its two delays are commensurate and its coefficients are constant. By contrast, the present characteristic equation contains two independently varying delays and the coefficient $e^{-\gamma\tau_2}$. The earlier integer-order studies \cite{hale1993global,belair1994stability,li1999stability,bhalekar2019analysing} provide important geometric and bifurcation results for multiple-delay DDEs, but they do not address the fractional power $\lambda^\alpha$ together with this delay-dependent gain.

Within this setting, the genuinely new results are the following. Theorem~\ref{thmcase1_unified_detailed} gives the stated delay-independent stability and instability regions, together with the critical quantities $\tau_1^*(\tau)$ and $\tau_{2*}$, for the reduced equation with the nonconstant coefficient $-k e^{-\gamma\tau}$. Theorem~\ref{thmcase2_unified} gives conditions that are uniform in both $\tau_1$ and $\tau_2$ and, in the indicated regions, uniform for $0<\alpha\leq1$. Theorem~\ref{Case2Thrm2} provides the explicit threshold $k_*$ that guarantees a positive characteristic root, and hence instability, in the fourth quadrant for $0<\alpha<1$. These conclusions cannot be obtained from the earlier constant-coefficient results by a direct substitution, because varying $\tau_2$ simultaneously changes a delay and a feedback coefficient.

The fractional order changes the spectral analysis in an essential way. For the integer-order case $\alpha=1$, the characteristic term is $\lambda$, and stability is governed by the usual left-half-plane condition. For $0<\alpha<1$, the term $\lambda^\alpha$ changes both the magnitude and the argument of each characteristic mode, and the corresponding sector condition depends explicitly on $\alpha$. In particular,
\[
(iv)^\alpha=v^\alpha\left[\cos\left(\frac{\alpha\pi}{2}\right)
+i\sin\left(\frac{\alpha\pi}{2}\right)\right],
\]
so a purely imaginary characteristic root contributes a nonzero real component when $\alpha<1$. The real and imaginary crossing equations, the critical-delay formulas, and the threshold $k_*$ therefore depend on $\alpha$ and do not reduce to the classical two-delay analysis. The distinction is especially visible in the fourth quadrant: when $\alpha=1$, $\lambda=-\gamma$ is automatically a positive characteristic root, whereas for $0<\alpha<1$ this identity is lost and instability must instead be established through the $\alpha$-dependent threshold in Theorem~\ref{Case2Thrm2}.

Recent work published in 2026 shows that well-posedness, memory effects, and multiple delays remain active themes in fractional dynamics. Abdalla, Liu, and Hammad \cite{abdalla2026existence} used fixed-point methods to establish existence results for impulsive Caputo delay systems. El-Saka, El-Sherbeny, and El-Sayed \cite{elsaka2026dynamic} studied fractional distributed-delay models and obtained stability regions from their characteristic equations. Zhu et al. \cite{zhu2026hopf} derived delay-induced Hopf conditions for a fractional-order Lotka--Volterra predator--prey model with two delays. These studies concern impulsive, distributed-delay, or multidimensional biological systems, whereas the present work treats a scalar Caputo equation with two discrete delays and the specific delay-dependent gain $e^{-\gamma\tau_2}$. The comparison highlights that the well-posedness result below, the explicit parameter regions, and the threshold $k_*$ address a structure not covered by those recent analyses.
}

This paper is organized as follows. Section~\ref{secprelim} reviews preliminary results on fractional delay differential equations with a single discrete delay. In Section~\ref{secmain}, we introduce the nonlinear fractional delay differential equation with delay-dependent coefficients, establish its equilibrium, and derive the corresponding linearized and characteristic equations. Sections~\ref{secC1} and~\ref{secC2} present the main stability theorems for the cases $\tau_1=0$ and $\tau_1>0$ and $\tau_2\ge0$ respectively, together with illustrative examples and graphical results. Finally, Section~\ref{secconclusion} provides concluding remarks and outlines possible directions for future research.

	\section{Preliminaries}
	\label{secprelim}
	We consider a linear delay differential equation with Caputo fractional derivative $D^\alpha$ \cite{podlubny1998fractional,kilbas2006theory}.
	\begin{theorem}\cite{bhalekar2016stability}
		\label{thmprelim}
		Consider the scalar FDE with a single discrete delay:
		\begin{equation}
			\label{eq1}
			D^\alpha x(t) = a x(t) + b x(t - \tau), \quad 0 < \alpha \le 1, \tau\ge0.
		\end{equation}
		Then, the zero equilibrium $x_*=0$ has the following stability behavior:
		\begin{enumerate}
			\item If $b \in (-\infty, -|a|)$, then the equilibrium is asymptotically stable for $\tau \in [0, \tau_{cr})$ and the system undergoes Hopf bifurcation at
			\[
			\tau_{cr} = \frac{\arccos\left( \frac{(a\cos(\frac{\alpha \pi}{2}) + \sqrt{b^2 - a^2 \sin^2(\frac{\alpha \pi}{2})})\cos(\frac{\alpha \pi}{2}) - a}{b} \right)}{\left( a \cos\left(\frac{\alpha \pi}{2}\right) + \sqrt{b^2 - a^2 \sin^2\left(\frac{\alpha \pi}{2}\right)} \right)^{1/\alpha}}.
			\]
			\item If $b \in (-a, \infty)$, then the equilibrium is unstable for all $\tau \ge 0$.
			\item If $a < 0$ and $b \in (a, -a)$, then the equilibrium is asymptotically stable for all $\tau \ge 0$.
		\end{enumerate}
		\begin{figure}[H]
			\centering
			\includegraphics[scale=0.57]{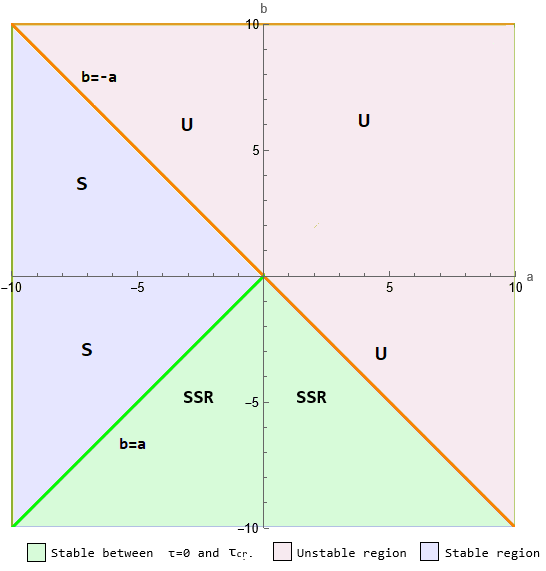}
			\caption{Stability regions of Equation (\ref{eq1}).}
			\label{figprel}
		\end{figure}
		\textnormal {Notation: Here, SSR means Single Stable Region.}
	\end{theorem}
	
	\section{The Main Model}
	\label{secmain}
	We consider the nonlinear fractional delay differential equation of the form:
	\begin{equation}
		D^\alpha x(t) = -\gamma x(t) + g(x(t - \tau_1)) - e^{-\gamma \tau_2} g(x(t - \tau_1 - \tau_2)), \quad 0 < \alpha \le 1,
		\label{eqmain}
	\end{equation}
	where $D^\alpha$ denotes the Caputo fractional derivative, $\gamma \in \mathbb{R}$, and $g \in C^1(\mathbb{R})$. The parameters $\tau_1, \tau_2 \ge 0$ represent discrete delays, and the exponential factor $e^{-\gamma \tau_2}$ introduces a delay-dependent coefficient. The integer-order counterpart of this equation, analyzed in \cite{boullu2020stability}, represents a model for platelet production.\\

	{
	\subsection*{Existence and uniqueness of solutions}
	Let $r=\tau_1+\tau_2$ and prescribe a history function
	\[
	x(t)=\phi(t),\qquad -r\leq t\leq0,
	\]
	where $\phi\in C([-r,0],\mathbb{R})$.

	\begin{theorem}\label{thmwellposed}
		Let $0<\alpha\leq1$, $\tau_1,\tau_2\geq0$, and $\gamma\in\mathbb{R}$. If $g:\mathbb{R}\to\mathbb{R}$ is locally Lipschitz continuous, then Equation~\eqref{eqmain}, together with the history $\phi$, has a unique solution on a maximal interval $[-r,T_{\max})$. If $g$ is globally Lipschitz continuous, then $T_{\max}=\infty$.
	\end{theorem}

	\begin{proof}
		For $t\geq0$, Equation~\eqref{eqmain} is equivalent to the Volterra integral equation
		\[
		x(t)=\phi(0)+\frac{1}{\Gamma(\alpha)}
		\int_0^t(t-s)^{\alpha-1}
		\left[
		-\gamma x(s)+g\bigl(x(s-\tau_1)\bigr)
		-e^{-\gamma\tau_2}g\bigl(x(s-r)\bigr)
		\right]\,ds.
		\]
		Fix a bounded neighborhood of the history and let $L$ be a Lipschitz constant of $g$ on that neighborhood. On
		\[
		X_h=\{x\in C([-r,h],\mathbb{R}):x(t)=\phi(t)\text{ for }-r\leq t\leq0\},
		\]
		define the operator $\mathcal{T}$ by the right-hand side of the integral equation. For $x,y\in X_h$,
		\[
		\|\mathcal{T}x-\mathcal{T}y\|_{\infty}
		\leq
		\frac{C h^\alpha}{\Gamma(\alpha+1)}
		\|x-y\|_{\infty},
		\qquad
		C=|\gamma|+L\bigl(1+e^{-\gamma\tau_2}\bigr).
		\]
		Choosing $h>0$ so that $C h^\alpha<\Gamma(\alpha+1)$ makes $\mathcal{T}$ a contraction on a sufficiently large closed ball in $X_h$. The Banach fixed-point theorem gives a unique solution on $[-r,h]$. Repeating the argument from the endpoint of each interval gives a unique maximal continuation on $[-r,T_{\max})$.

		If $g$ is globally Lipschitz with constant $L$, then
		\[
		|g(u)|\leq L|u|+|g(0)|
		\]
		for all $u\in\mathbb{R}$. Consequently, on every finite interval $[0,T]$, the function
		\[
		M(t)=\sup_{-r\leq s\leq t}|x(s)|
		\]
		satisfies an inequality of the form
		\[
		M(t)\leq A_T+B\int_0^t(t-s)^{\alpha-1}M(s)\,ds,
		\]
		where $A_T,B>0$. The generalized fractional Gronwall inequality of Ye, Gao, and Ding \cite{ye2007generalized} then yields a Mittag--Leffler bound for $M(t)$. A retarded version adapted specifically to fractional differential equations with delay was established by Ye and Gao \cite{ye2011henry}; see also Diethelm \cite{diethelm2010analysis} for the corresponding Caputo initial-value theory. Hence, $|x(t)|$ remains bounded on every finite interval. Finite-time blow-up cannot occur, and the continuation extends to all $t\geq0$. This proves global existence and uniqueness.
	\end{proof}

	\subsection*{Positive-parameter physical interpretation}
	For the biologically relevant choice
	\[
	\gamma>0,\qquad \tau_1>0,\qquad \tau_2>0,
	\]
	the parameters have a direct interpretation. In the platelet-production model of Boullu, Pujo-Menjouet, and B{\'e}lair \cite{boullu2020stability}, $x(t)$ represents a mature cell or platelet population, $\gamma$ is a positive removal rate, and the two positive delays describe successive biological stages. The factor $e^{-\gamma\tau_2}$ is the survival fraction accumulated during the second delayed stage. Consequently, the exact integer-order model underlying Equation~\eqref{eqmain} already occurs in a calibrated physiological application; varying $\tau_2$ or $\gamma$ can move the platelet equilibrium across a stability boundary. Related positive-parameter two-delay structures also arise in age-structured erythropoiesis \cite{belair1995age}, solid-tumor growth \cite{piotrowska2008hopf}, and, in recent fractional modeling, predator--prey interactions with two positive response delays \cite{zhu2026hopf}. These examples support the use of positive decay, interaction, and delay parameters in the analytical regions studied below.
	}

	Let $x_*$ be an equilibrium of Equation~\eqref{eqmain}. It must satisfy
	\[	-\gamma x_* + g( x_*) - g( x_*) e^{-\gamma \tau_2} = 0. \]
	We assume that $g$ is nonlinear and satisfies $g(0)=0$. Under this assumption, $x_*=0$ is an equilibrium of Equation~\eqref{eqmain}.
	Hence, assuming $ g'(0)=k$, we linearize equation (\ref{eqmain}) about the equilibrium $x_*= 0$.
	\begin{enumerate}

		\item \textbf{Linearization near Equilibrium}

		Let $x(t)$ be an infinitesimally perturbed solution of (\ref{eqmain}) near $x_*=0.$
		Using the first-order Taylor approximation for the nonlinear terms about the equilibrium and denoting $g'(0) = k$, we have
		\[g(x(t - \tau_1)) \approx g(0) + g'(0) x(t - \tau_1) = k x(t - \tau_1),\]
		\[g(x(t - \tau_1 - \tau_2)) \approx g(0) + g'(0) x(t - \tau_1 - \tau_2) = k x(t - \tau_1 - \tau_2).\]
		
		Substituting into equation~\eqref{eqmain}, we obtain the \textbf{linearized equation} i.e.
		\begin{equation}
			D^\alpha x(t) = -\gamma x(t) + k x(t - \tau_1) - k e^{-\gamma \tau_2} x(t - \tau_1 - \tau_2).
			\label{eqlinear}
		\end{equation}
		We analyze the stability of this linearized equation in the presiding sections for various particular cases of parameters.
		
		\item \textbf {Characteristic Equation}
		
		To analyze stability, we apply the Laplace transform and obtain the characteristic equation corresponding to equation~\eqref{eqlinear} as
		\begin{equation}
			\lambda^\alpha = -\gamma + k e^{-\lambda \tau_1} - k e^{-\gamma \tau_2} e^{-\lambda (\tau_1 + \tau_2)}.
			\label{eqchar}
		\end{equation}

		{
		The fractional power in Equation~\eqref{eqchar} is not a formal replacement of $\lambda$ by $\lambda^\alpha$. It changes the stability boundary because, for $0<\alpha<1$, the argument of a characteristic mode is multiplied by $\alpha$. Accordingly, the non-delayed fractional equation is asymptotically stable when its characteristic values satisfy the sector condition
		\[
		|\arg(\lambda)|>\frac{\alpha\pi}{2},
		\]
		rather than only the integer-order condition $\operatorname{Re}(\lambda)<0$. At a possible imaginary-axis crossing, the identity
		\[
		(iv)^\alpha
		=v^\alpha\cos\left(\frac{\alpha\pi}{2}\right)
		+i\,v^\alpha\sin\left(\frac{\alpha\pi}{2}\right)
		\]
		produces the explicit $\alpha$-dependent terms used in Equations~\eqref{eqreal_case2} and~\eqref{eqimag_case2}. Thus, changing $\alpha$ changes the crossing frequency, the critical delays, and the feedback threshold. When $\alpha=1$, the cosine term vanishes and the classical integer-order crossing equations are recovered.
		}
		
	\end{enumerate}

	\section{Case 1: $\tau_1 = 0$}
	\label{secC1}
	Setting $\tau_1=0$ and $\tau_2=\tau$ in Equation~\eqref{eqlinear} yields
	\begin{equation}
		D^\alpha x(t) = (k - \gamma) x(t) - k e^{-\gamma \tau} x(t - \tau), \quad 0 < \alpha \le 1.
		\label{eqcase1}
	\end{equation}
    Note that even though the Equation (\ref{eqcase1}) is a scalar DDE with a single delay, it is not same as that of Equation (\ref{eq1}) studied in the literature. The coefficient of $x(t-\tau)$ in (\ref{eqcase1}) depends on the delay, unlike (\ref{eq1}). The Equation (\ref{eqcase1}) is of the general form 
	\begin{equation}
		\label{c1lin}
		D^\alpha x(t) = a x(t) + b(\tau) x(t - \tau),
	\end{equation}
	where $a = k - \gamma$ and $b(\tau) = -k e^{-\gamma \tau}$ depends on the delay.
	We utilize the existing stability results (Theorem \ref{thmprelim}) for the general form (\ref{c1lin}) for any fixed value of $\tau$ to determine regions of delay-independent stability and instability.

	The following theorem summarizes the stability regions of the equilibrium $x_* = 0$ of (\ref{eqmain}) for this case in the $k-\gamma$ plane. 
	\begin{theorem}
		\label{thmcase1_unified_detailed}
		Consider the equation~\eqref{eqcase1} with parameters $k, \gamma \in \mathbb{R}$ and $\tau \ge 0$. Then the equilibrium $x_* = 0$ has the following behaviors (cf. Figure ~\ref{figcase1}):
		\begin{enumerate}[label=(\alph*)]
			
			\item \textbf{Delay-Independent instability.}  
			If either of the following holds:
			\begin{enumerate}[label=(\roman*)]
				\item $\gamma < k < 0$, or
				\item $k < \gamma < 0$,
			\end{enumerate}
			i.e., $(k,\gamma)$ lies in the third quadrant, then the equilibrium is unstable for all $\tau \ge 0$.
			\item \textbf{Delay-Independent stability.}  
			The equilibrium is asymptotically stable for all $\tau \ge 0$ if either:
			\begin{enumerate}[label=(\roman*)]
				\item $\gamma > 2k > 0$, or
				\item $\gamma > 0$ and $k < 0$ i.e. the second quadrant of $k-\gamma$ plane.
			\end{enumerate}
			\item \textbf{Delay-Dependent stability.} 
			\begin{enumerate}[label=(\roman*)]
				\item Let \( 0 < \gamma < 2k \).
				Define \[\tau_{2*}= 
				\frac{-1}{\gamma} \log\left( \frac{|k - \gamma|}{k} \right),\;\;k\ne\gamma  \;\;\text{and}\]   \[\tau_1^*(\tau) = \frac{\arccos\left( \dfrac{(k-\gamma)\cos\left(\frac{\alpha \pi}{2}\right) + \sqrt{(-ke^{-\gamma\tau})^2 - (k-\gamma)^2 \sin^2\left(\frac{\alpha \pi}{2}\right)} \cdot \cos\left(\frac{\alpha \pi}{2}\right) - (k-\gamma)}{-ke^{-\gamma\tau}} \right)}{\left( (k-\gamma)\cos\left(\frac{\alpha \pi}{2}\right) + \sqrt{(-ke^{-\gamma\tau})^2 -(k-\gamma)^2 \sin^2\left(\frac{\alpha \pi}{2}\right)} \right)^{1/\alpha}}\]
				Then the stability of the equilibrium is given by:
				\begin{itemize}
					\item For \( \tau < \tau_2^* \), the system is:
					\begin{itemize}
						\item \textbf{Stable} if \( \tau < \tau_1^*(\tau) \),
						\item \textbf{Unstable} if \( \tau > \tau_1^*(\tau) \).
					\end{itemize}
					\item For \( \tau > \tau_2^* \), the system is:
					\begin{itemize}
						\item \textbf{Unstable} if \( 0 < \gamma < k \),
						\item \textbf{Stable} if \( k < \gamma < 2k \).
					\end{itemize}

				\end{itemize}
				
				\item Let $k>0,\gamma<0\; \text{i.e. }$ fourth quadrant of the $k-\gamma$ plane. Then, the stability is given by:
				\begin{itemize}
					\item For $\tau<\tau_2^*$, equilibrium point is unstable
					\item For $\tau>\tau_2^*$, the equilibrium is 
					\begin{itemize}
						\item stable if $\tau<\tau_1^*(\tau)$,
						\item unstable if $\tau>\tau_1^*(\tau)$.
					\end{itemize}
				\end{itemize}
			\end{enumerate}
		\end{enumerate}
	\end{theorem}
	
	\begin{proof}
		\textbf{(a)}
		In both subcases (a)(i) and (a)(ii), $k<0$ and $\gamma<0$. Hence, for every $\tau>0$, one has $e^{-\gamma\tau}>1$ and
		\[
		b(\tau)=-k e^{-\gamma\tau}>0.
		\]
		
		\textbf{Case (i):} If $\gamma<k<0$, then $a=k-\gamma>0$. Consequently, $-a<0<b(\tau)$, and hence
		
		\[
		b(\tau)\in(-a,\infty).
		\]
		Theorem~\ref{thmprelim} therefore implies that the equilibrium is unstable for every $\tau\geq0$.\\
		
		\textbf{Case (ii):} If $k<\gamma<0$, then $a=k-\gamma<0$. Define
		\[
		f(\tau):=b(\tau)+a=-k e^{-\gamma\tau}+k-\gamma.
		\]
		At $\tau=0$,
		\[
		f(0) = -k + k - \gamma = -\gamma > 0.
		\]
		Moreover,
		\[
		f'(\tau) = k \gamma e^{-\gamma \tau} > 0 \quad (\text{since } k< \gamma < 0).
		\]
		Thus, $f$ is strictly increasing on $[0,\infty)$, and $f(\tau)\geq f(0)>0$ for every $\tau\geq0$. It follows that
		\[
		b(\tau) > -a.
		\]
		Therefore, $b(\tau)\in(-a,\infty)$, and Theorem~\ref{thmprelim} shows that the equilibrium is unstable for every $\tau\geq0$.\\
		
		\textbf{(b)}
		We show that $a < 0$ and $|b(\tau)| < -a$ in both cases b(i) and b(ii), which will ensure the global stability via Part (3) of Theorem~\ref{thmprelim}.
		
		\textbf{Case (i):} Suppose that $\gamma>2k>0$. Then $k>0$ and $a=k-\gamma<0$. Since $\gamma>0$, for every $\tau>0$,
		
		\[
		|b(\tau)|=k e^{-\gamma\tau}<k.
		\]
		
		The assumption $\gamma>2k$ also gives $-a=\gamma-k>k$. Hence,
		\[
		|b(\tau)| < -a \quad \text{for all } \tau>0.
		\]
		Part (3) of Theorem~\ref{thmprelim} now shows that $x_*=0$ is asymptotically stable for every $\tau\geq0$.\\
		
		\textbf{Case (ii):} Suppose that $\gamma>0$ and $k<0$. Then $a=k-\gamma<0$ and $b(\tau)=-k e^{-\gamma\tau}>0$. Furthermore,
		
		Thus,
		\[
		|b(\tau)| = -k e^{-\gamma \tau} \le -k < \gamma - k = -a.
		\]
		Thus, $|b(\tau)|<-a$, so Part (3) of Theorem~\ref{thmprelim} proves asymptotic stability for every $\tau\geq0$.\\
		
		\textbf{(c)}\\ 
		\textbf{Case (i):} Assume that $0<\gamma<2k$. Then $2k-\gamma>0$, and therefore $k-\gamma>-k$. Since $a=k-\gamma$ and $b(0)=-k$, this inequality is equivalent to
		\begin{equation}
			a>b(0). \tag{*}\label{star}
		\end{equation}
		At $\tau = 0$, equation (\ref{eqcase1}) becomes
		\[
		D^\alpha x(t) = -\gamma x(t).
		\]
		Since $\gamma > 0$, equilibrium is asymptotically stable at $\tau = 0$.
		We next verify the position of $(a,b(0))$ in the stability diagram.
		\begin{itemize}
			\item If $a<0$, then $a=-|a|$. Inequality~\eqref{star} gives $b(0)<a=-|a|$.
			\item If $a=k-\gamma>0$, then $k>\gamma>0$. Hence $-k<-k+\gamma=-a$, so $b(0)<-a=-|a|$.
		\end{itemize}
		In either case, $b(0)<-|a|$. Thus, $(a,b(0))$ lies in the single stable region shown in Figure~\ref{figprel}.
		\vspace{0.3cm}
		The function $b(\tau)=-k e^{-\gamma\tau}$ is strictly increasing and converges to $0$ as $\tau\to\infty$. Hence, for sufficiently small $\tau$, the point $(a,b(\tau))$ remains in the single stable region of Figure~\ref{figprel}. Applying Part (1) of Theorem~\ref{thmprelim} with $a=k-\gamma$ and $b=-k e^{-\gamma\tau}$ gives the critical value $\tau_1^*(\tau)$ displayed in the theorem statement. The second critical value, $\tau_{2*}$, is the value at which $(a,b(\tau))$ reaches the relevant boundary of the stability region; see Figure~\ref{case1}.
		\begin{figure}[H]
			\centering
			\includegraphics[scale=0.55]{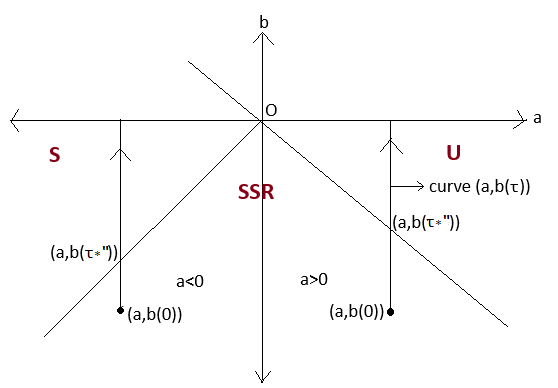} 
			\caption{Analysis for the region $0<\gamma<2k$ }
			\label{case1}
		\end{figure}
		For $\tau<\tau_{2*}$, Part (1) of Theorem~\ref{thmprelim} gives the following alternatives:
		\begin{itemize}
			\item equilibrium is stable if $\tau<\tau_1^*(\tau)$ and 
			\item unstable if $\tau>\tau_1^*(\tau)$ (by Theorem \ref{thmprelim}, case 1).
		\end{itemize}
		If $a>0$, the boundary is $b=-a$. Solving
		\[
		-k e^{-\gamma\tau_{2*}}=-a=\gamma-k
		\]
		for $\tau_{2*}$ yields
		\[
		\tau_{2*}=-\frac{1}{\gamma}\log\left(\frac{k-\gamma}{k}\right).
		\]
		If $a<0$, the boundary is $b=a$. Solving $-k e^{-\gamma\tau_{2*}}=k-\gamma$ gives
		\[
		\tau_{2*}=-\frac{1}{\gamma}\log\left(\frac{\gamma-k}{k}\right).
		\]
		The two formulas can be written as
		\[
		\tau_{2*}=-\frac{1}{\gamma}\log\left(\frac{|k-\gamma|}{k}\right),\qquad k\ne\gamma.
		\]
		\begin{figure}[H]
			\centering
			\includegraphics[scale=0.45]{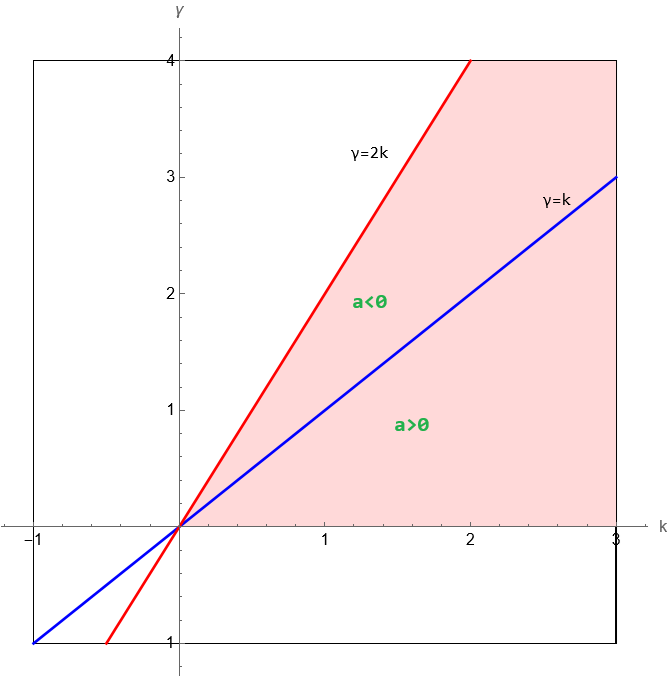}
			\caption{Region: $0<\gamma<2k$}
			\label{case1dd}
		\end{figure}
		For $\tau>\tau_{2*}$, Figure~\ref{case1} and Theorem~\ref{thmprelim} show that the equilibrium is stable when $a<0$ and unstable when $a>0$. In the present parameter region, $a<0$ is equivalent to $k<\gamma<2k$, whereas $a>0$ is equivalent to $0<\gamma<k$. Therefore, for $\tau>\tau_{2*}$, the equilibrium is asymptotically stable if $k<\gamma<2k$ and unstable if $0<\gamma<k$; see Figure~\ref{case1dd}.\\
		
		\textbf{Case (ii):} $k > 0$ and $\gamma < 0$ \\
		At $\tau = 0$, equation~\eqref{eqmain} reduces to
		\[
		D^\alpha x(t) = -\gamma x(t).
		\]
		Since $\gamma<0$, the equilibrium is unstable at $\tau=0$. Moreover, $a=k-\gamma>0$ and $b(\tau)=-k e^{-\gamma\tau}<0$. The inequality $-a=\gamma-k<-k=b(0)$ shows that $(a,b(0))$ lies in the unstable region of Figure~\ref{figprel}.
		
		Because $b'(\tau)=k\gamma e^{-\gamma\tau}<0$, the function $b$ is strictly decreasing. It crosses the boundary $b=-a$ at the unique value
		\[
		\tau_2^* = \frac{-1}{\gamma} \log\left( \frac{k-\gamma}{k} \right).
		\]
		For $\tau<\tau_2^*$, the point $(a,b(\tau))$ remains in the unstable region. For $\tau>\tau_2^*$, it lies in the single stable region in which stability is determined by the critical-delay condition in Part (1) of Theorem~\ref{thmprelim}.
		
		Substitution of $a=k-\gamma$ and $b=-k e^{-\gamma\tau}$ into the expression for $\tau_{cr}$ in Theorem~\ref{thmprelim} yields the stated function $\tau_1^*(\tau)$.
		
		Consequently, for $\tau>\tau_2^*$, the equilibrium is stable when $\tau<\tau_1^*(\tau)$ and unstable when $\tau>\tau_1^*(\tau)$. Thus, whenever these inequalities are satisfied successively as $\tau$ increases, the system exhibits an unstable--stable--unstable sequence.
	\end{proof}

	{
	\subsection*{Connection with the Hopf condition}
	For $\tau_1=0$, define
	\[
	\Delta(\lambda,\tau)
	=\lambda^\alpha-(k-\gamma)+k e^{-(\gamma+\lambda)\tau}.
	\]
	In the delay-dependent region $b(\tau)=-k e^{-\gamma\tau}<-|k-\gamma|$, Theorem~\ref{thmprelim} applies at each fixed value of the coefficient. The boundary
	\[
	\tau=\tau_1^*(\tau)
	\]
	is precisely the condition under which $\Delta(iv,\tau)=0$ for some $v>0$, and hence it is the Hopf candidate inherited from the critical delay $\tau_{cr}$ in Theorem~\ref{thmprelim}. If the root pair $\lambda=\pm iv$ is simple and crosses the imaginary axis transversally, then the stability change at this boundary is a Hopf bifurcation. The transversality condition can be checked from
	\[
	\frac{d\lambda}{d\tau}
	=
	\frac{k(\gamma+\lambda)e^{-(\gamma+\lambda)\tau}}
	{\alpha\lambda^{\alpha-1}-k\tau e^{-(\gamma+\lambda)\tau}},
	\]
	by verifying that
	\[
	\operatorname{Re}\left(\frac{d\lambda}{d\tau}\right)_{\lambda=iv,\ \tau=\tau_c}\neq0.
	\]
	Thus, the inequalities $\tau<\tau_1^*(\tau)$ and $\tau>\tau_1^*(\tau)$ describe the stable and unstable sides of the Hopf boundary, respectively, whenever the simplicity and transversality conditions hold. By contrast, $\tau_{2*}$ marks the crossing of the coefficient curve $b(\tau)$ through a boundary of the $(a,b)$ stability diagram; it is not, by itself, a Hopf condition. In the delay-independent stable and unstable regions of Theorem~\ref{thmcase1_unified_detailed}, no Hopf crossing occurs as $\tau$ varies.
	}

	\begin{figure}[H]
		\centering
		\includegraphics[scale=0.5]{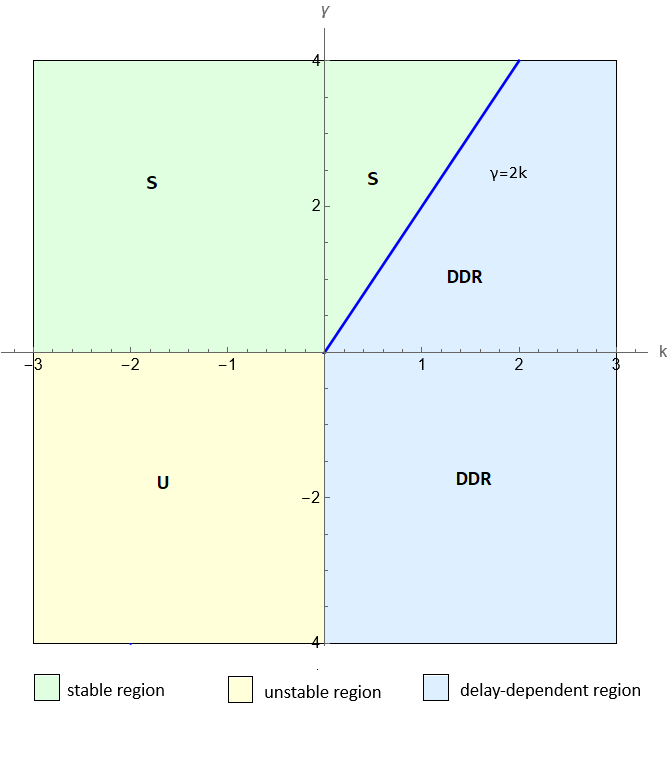}
		\caption{Stability regions for Case 1}
		\label{figcase1}
	\end{figure}
	\subsection{Illustrative Examples for Case 1}
	Consider FDDE $D^{\alpha} x(t) = (k-\gamma) x(t) + k e^{-\gamma\tau} x(t - \tau).$\\
	
\begin{figure}[H]
		\centering
		\includegraphics[scale=0.4]{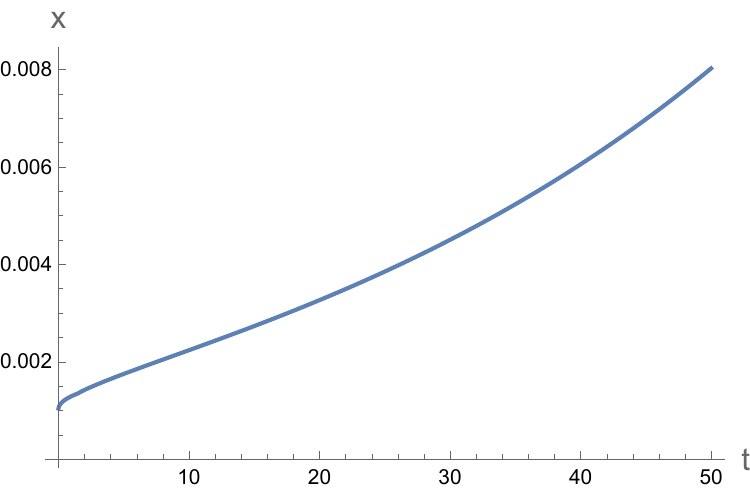}
		\caption{$\alpha=0.4,\; k=-0.1,\; \gamma=-0.2,\; \tau=1.5$}
		\label{1a}
	\end{figure}

	\textbf{Example 4.1}
	Suppose $\alpha=0.4,\; k=-0.1,\; \gamma=-0.2,\; \tau=1.5$.  
	Here, $k,\gamma<0$. By Theorem~\ref{thmcase1_unified_detailed}(a), the equilibrium is delay-independently unstable (cf. Figure \ref{1a}).
	
\begin{figure}[H]
		\centering
		\includegraphics[scale=0.4]{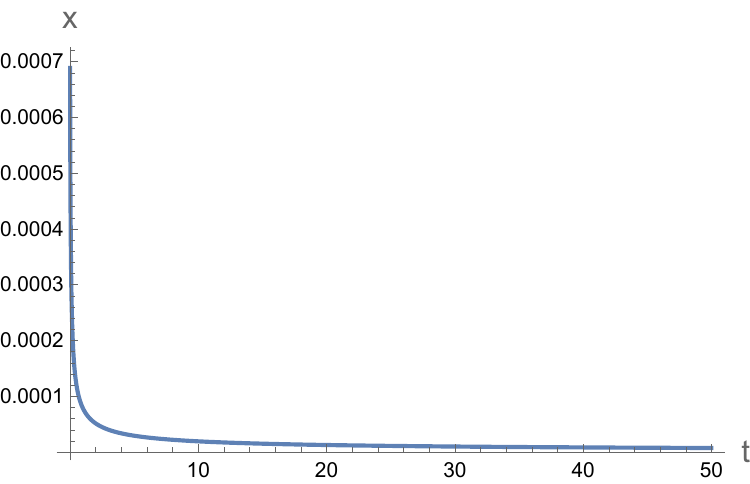}
		\caption{$\alpha=0.6,\; k=4,\; \gamma=10,\; \tau=0.8$}
		\label{1b}
	\end{figure}

	\textbf{Example 4.2}
	Now, consider $\alpha=0.6,\; k=4,\; \gamma=10,\; \tau=0.8$.  
	Here $\gamma>2k>0$. By Theorem~\ref{thmcase1_unified_detailed}(b)(i), the equilibrium is delay-independently stable. Hence, Numerical solutions converge to zero (see Figure \ref{1b}).
	
\begin{figure}[H]
		\centering
		\includegraphics[scale=0.4]{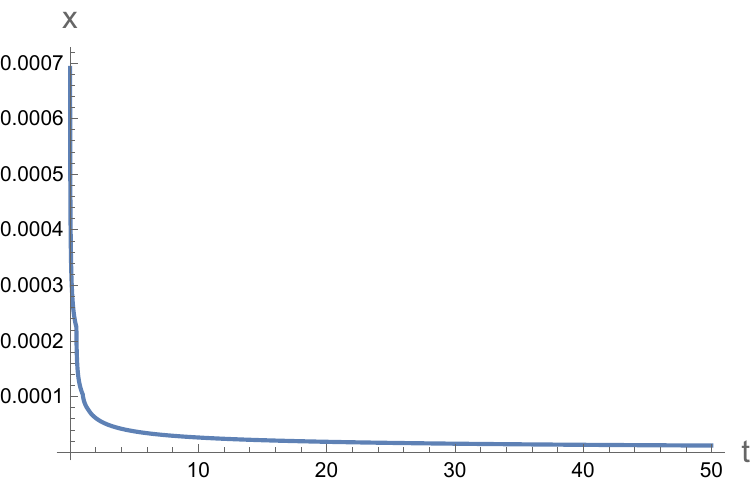}
		\caption{ $\alpha=0.5,\; k=-5,\; \gamma=3,\; \tau=0.5$}
		\label{1c}
	\end{figure}

	\textbf{Example 4.3}
	Take $\alpha=0.5,\; k=-5,\; \gamma=3,\; \tau=0.5$.  
	Here $k<0,\;\gamma>0$. Theorem~\ref{thmcase1_unified_detailed}(b)(ii) ensures global stability, which is verified by numerical simulations (cf. Figure \ref{1c}).
	
\begin{figure}[H]
		\centering
		\includegraphics[scale=0.4]{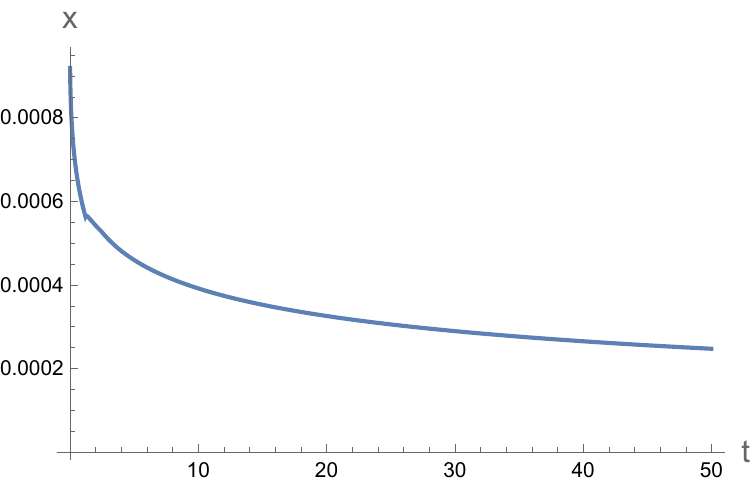}
		\caption{$\alpha=0.4,\; k=1.5,\; \gamma=1.8,\; \tau=1.2$}
		\label{1d}
	\end{figure}

	\textbf{Example 4.4}
	For $\alpha=0.4,\; k=1.5,\; \gamma=1.8,\; \tau=1.2$.  
	Here, $0<\gamma<2k$ and $\tau_{2*}=0.8941$. Since, $\tau>\tau_{2*}$, Theorem~\ref{thmcase1_unified_detailed}(c) predicts stability (see Figure \ref{1d}). \\
	Now, suppose $\tau=0.3<\tau_{2*}$, so we find $\tau_{1*}(\tau)=9.2161$. Again, $\tau<\tau_{1*}$, so the equilibrium is stable.
	\begin{figure}[H]
		\centering
		\includegraphics[scale=0.4]{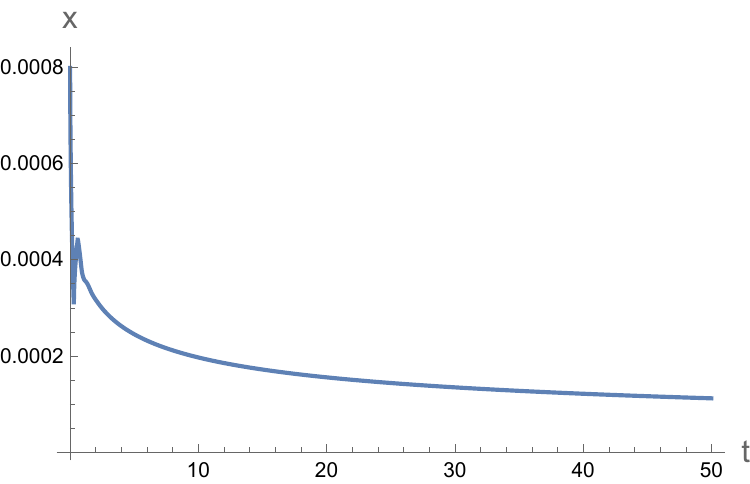}
		\caption{$\alpha=0.4,\; k=1.5,\; \gamma=1.8,\; \tau=0.3$}
	\end{figure}
	
\begin{figure}[H]
		\centering
		\includegraphics[scale=0.4]{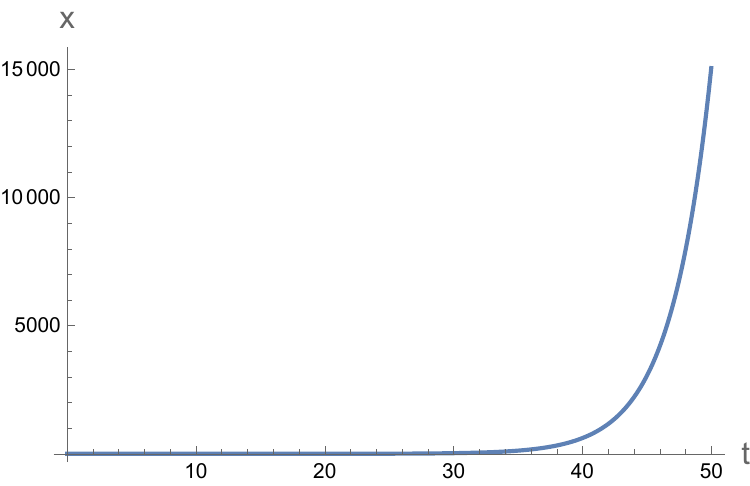}
		\caption{ $\alpha=0.7,\; k=3.4,\; \gamma=2.2,\; \tau=0.6$}
		\label{1f}
	\end{figure}

	\textbf{Example 4.5}
	Take, $\alpha=0.7,\; k=3.4,\; \gamma=2.2,\; \tau=0.6$.  
	Here, again $0<\gamma<2k$ and $\tau_{2*}=0.4733 $. Here, $\tau>\tau_{2*}$, which implies instability. Numerical Simulations verify this (cf. Figure (\ref{1f}).
	Now, let $\tau=0.2<\tau_{2*}$, so we find $\tau_{1*}(\tau)=0.4244$. Again, $\tau<\tau_{1*}$, so the equilibrium is stable.
	\begin{figure}[H]
		\centering
		\includegraphics[scale=0.4]{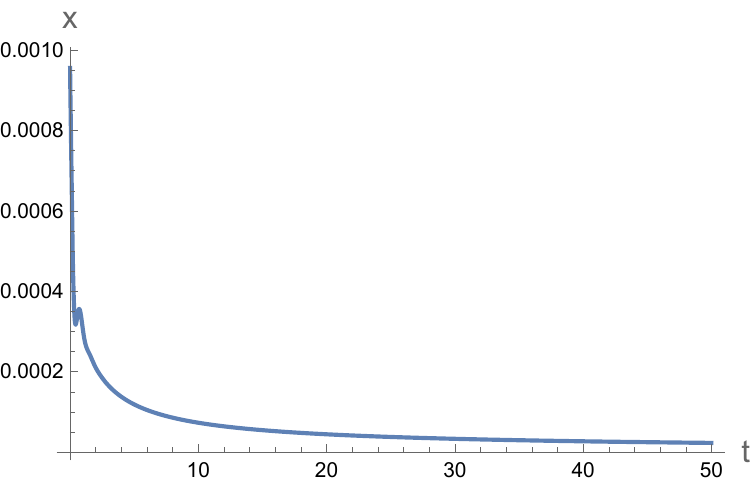}
		\caption{ $\alpha=0.7,\; k=3.4,\; \gamma=2.2,\; \tau=0.2$}
	\end{figure}
	\textbf{Example 4.6}
	For $\alpha=0.8,\; k=3,\; \gamma=-5,\; \tau=0.08$.  
	Here $k>0$ and $\gamma<0$, and the critical value is $\tau_2^*=0.1962$. Since $\tau<\tau_2^*$, the theorem predicts instability, which is confirmed numerically in Figure~\ref{1h}.
	\begin{figure}[H]
		\centering
		\includegraphics[scale=0.4]{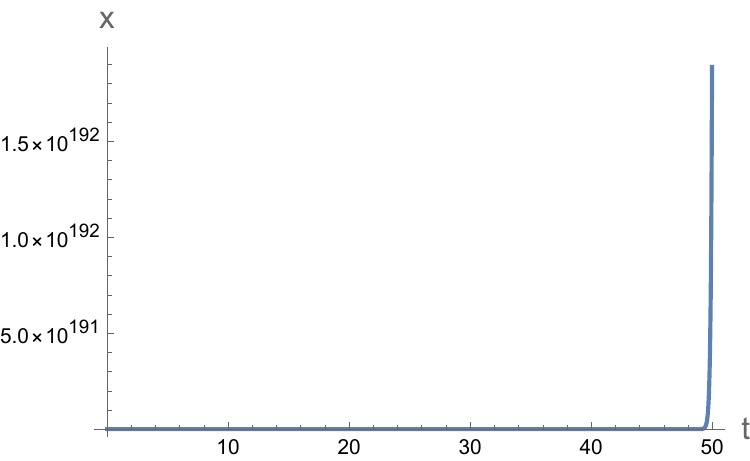}
		\caption{$\alpha=0.8,\; k=3,\; \gamma=-5,\; \tau=0.08$}
		\label{1h}
	\end{figure}
	Now, let  $\tau=0.35>\tau_2^*$;  $\tau_{1*}(\tau)=0.038614$. Since, $\tau>\tau_{1*}$ implies an unstable equilibrium point (Figure \ref{1i}).
	\begin{figure}[H]
		\centering
		\includegraphics[scale=0.4]{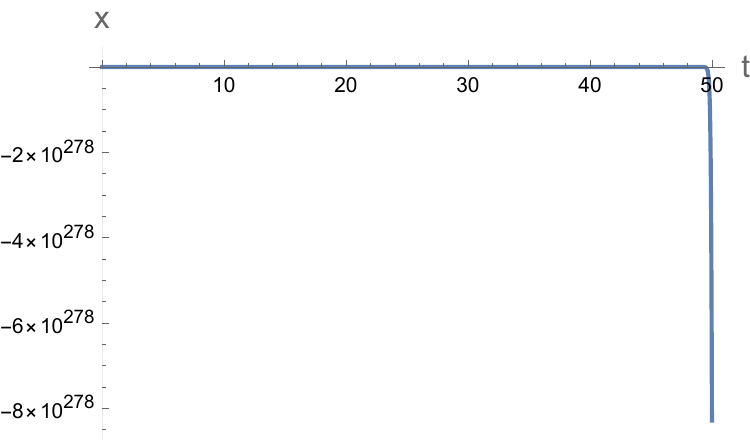}
		\caption{$\alpha=0.8,\; k=3,\; \gamma=-5,\; \tau=0.35$}
		\label{1i}
	\end{figure}
	\textbf{Note:} By testing a range of $\tau$ values, we verify that stability is independent of delay; however, for the sake of conciseness, we only show a few representative examples here.\\

	\section{Case 2: $\tau_1>0, \;\tau_2\ge0$}
	\label{secC2}
	In this case, both delays are retained. The linearized version of this equation is given by (\ref{eqlinear}).
	The corresponding characteristic equation is given by $(\ref{eqchar})$.

	We now present some delay-independent results summarizing stability behaviors for this case.
	\begin{theorem}
		\label{thmcase2_unified}
		The equilibrium point $x_* = 0$ of Equation (\ref{eqmain}) exhibits the following stability behaviors (see Fig. \ref{ind2}):
		\begin{enumerate}[label=(\alph*)]
			\item \textbf{Delay-Independent Stability:} The equilibrium is asymptotically stable for all delay values $\tau_1\ge0,\;\tau_2\ge0$ and for $0<\alpha\le1$. if either
			\begin{enumerate}[label=(\roman*)]
				\item $\gamma > 2k > 0$, or
				\item $\gamma > -2k > 0$ .
			\end{enumerate}
			\item \textbf{Delay-Independent Instability:}
			\begin{enumerate}[label=(\roman*)]
				\item If $\gamma < 0$ and $k < 0$, i.e. $(k,\gamma)$ lies in the third quadrant of the $k\gamma$-plane, then the equilibrium is unstable for all $\tau_1\ge0,\;\tau_2\ge0$ and for $0<\alpha\le1$.
				\item For $\alpha=1,$ the equilibrium is unstable in the fourth quadrant of the $k\gamma$-plane.
			\end{enumerate}
		\end{enumerate}
	\end{theorem}
	\begin{proof}
		\textbf{(a)} Fix $\tau_1\geq0$ and use the principal branch of $\lambda^\alpha$. When $\tau_2=0$, the two delayed terms cancel, and Equation~\eqref{eqlinear} reduces to
		\[
		D^\alpha x(t)=-\gamma x(t).
		\]
		Its characteristic equation is $\lambda^\alpha=-\gamma$. Since $\gamma>0$ in both cases (i) and (ii), the value $-\gamma$ has argument $\pi$, which satisfies
		\[
		\pi>\frac{\alpha\pi}{2},\qquad 0<\alpha\leq1.
		\]
		Thus, the equilibrium is asymptotically stable at $\tau_2=0$.

		We next exclude every root on the imaginary axis for $\tau_2>0$. First consider the zero-frequency case. Substitution of $\lambda=0$ into Equation~\eqref{eqchar} gives
		\[
		0=-\gamma+k\left(1-e^{-\gamma\tau_2}\right).
		\]
		In case (i), $k>0$ and $0<1-e^{-\gamma\tau_2}<1$, so
		\[
		-\gamma+k\left(1-e^{-\gamma\tau_2}\right)
		<-\gamma+k<0.
		\]
		In case (ii), $k<0$ and $1-e^{-\gamma\tau_2}>0$, so
		\[
		-\gamma+k\left(1-e^{-\gamma\tau_2}\right)<-\gamma<0.
		\]
		Hence, $\lambda=0$ is not a characteristic root in either case.

		It remains to exclude roots $\lambda=iv$ with $v\neq0$. Since the characteristic equation has real coefficients, nonreal roots occur in conjugate pairs; it is therefore sufficient to consider $v>0$. Substitution of $\lambda=iv$ into Equation~\eqref{eqchar} gives
		\[
		(iv)^\alpha
		=-\gamma+k e^{-iv\tau_1}
		-k e^{-\gamma\tau_2}e^{-iv(\tau_1+\tau_2)}.
		\]
		Using
		\[
		(iv)^\alpha
		=v^\alpha\left[
		\cos\left(\frac{\alpha\pi}{2}\right)
		+i\sin\left(\frac{\alpha\pi}{2}\right)
		\right]
		\]
		and Euler's formula for the exponential terms yields
		\begin{align}
			v^\alpha \cos\left( \frac{\alpha \pi}{2} \right)
			&= -\gamma + k \cos(v \tau_1)
			- k e^{-\gamma \tau_2}
			\cos\bigl(v(\tau_1+\tau_2)\bigr), \label{eqreal_case2}\\
			v^\alpha \sin\left( \frac{\alpha \pi}{2} \right)
			&= -k \sin(v \tau_1)
			+ k e^{-\gamma \tau_2}
			\sin\bigl(v(\tau_1+\tau_2)\bigr). \label{eqimag_case2}
		\end{align}
		Any nonzero imaginary root must satisfy both equations. We show that the real equation alone is impossible in the parameter regions of part (a).

		\textbf{Case (i):} Suppose that $\gamma>2k>0$. Put
		\[
		E=e^{-\gamma\tau_2},\qquad
		\theta_1=v\tau_1,\qquad
		\theta_2=v(\tau_1+\tau_2).
		\]
		Because $\gamma>0$ and $\tau_2>0$, one has $0<E<1$. Since $k>0$ and $-1\leq\cos\theta\leq1$,
		\[
		k\cos\theta_1\leq k
		\]
		and
		\[
		-kE\cos\theta_2
		\leq kE.
		\]
		The second inequality follows from $-\cos\theta_2\leq1$ and the positivity of $kE$. Adding the two inequalities gives
		\[
		k\cos\theta_1-kE\cos\theta_2
		\leq k+kE<2k.
		\]
		Consequently, the right-hand side of Equation~\eqref{eqreal_case2} satisfies
		\[
		-\gamma+k\cos\theta_1-kE\cos\theta_2
		<-\gamma+2k<0.
		\]
		By contrast,
		\[
		v^\alpha\cos\left(\frac{\alpha\pi}{2}\right)\geq0
		\]
		for $v>0$ and $0<\alpha\leq1$. The two sides of Equation~\eqref{eqreal_case2} therefore cannot be equal, so no nonzero purely imaginary root exists.

		\textbf{Case (ii):} Suppose that $\gamma>-2k>0$, so $k<0$. Again let $E$, $\theta_1$, and $\theta_2$ be as above. Since $\cos\theta_1\geq-1$, multiplication by the negative number $k$ reverses the inequality and gives
		\[
		k\cos\theta_1\leq-k.
		\]
		Furthermore, $-kE>0$ and $\cos\theta_2\leq1$, so
		\[
		-kE\cos\theta_2\leq-kE.
		\]
		Hence,
		\[
		k\cos\theta_1-kE\cos\theta_2
		\leq-k-kE.
		\]
		Because $0<E<1$ and $-k>0$, one has $-kE<-k$, and therefore
		\[
		-\gamma-k-kE<-\gamma-2k<0.
		\]
		It follows that the right-hand side of Equation~\eqref{eqreal_case2} is strictly negative, while its left-hand side is nonnegative. Thus, no nonzero purely imaginary root exists in case (ii).

		Finally, connect $\tau_2=0$ to an arbitrary target value $\tau_2>0$ through the homotopy $s\mapsto s\tau_2$, $0\leq s\leq1$, while keeping $\tau_1$ fixed. By the standard root-continuation principle for retarded characteristic equations, characteristic roots depend continuously on the delay parameter on bounded subsets of the complex plane, and the number of roots in the open right half-plane can change only when a root reaches the imaginary axis \cite{hale2013introduction,bhalekar2016stability}. Since the equilibrium is stable at $s=0$, and the preceding argument excludes both zero and nonzero imaginary roots for every $s\in(0,1]$, no root can enter the closed right half-plane along this path. Therefore, the equilibrium remains asymptotically stable for all $\tau_1,\tau_2\geq0$ in both parameter regions.\\
		
		\textbf{(b)} \\
		\textbf{Case (i)} Assume that $\gamma<0$ and $k<0$. At $\tau_2=0$, the characteristic equation is $\lambda^\alpha=-\gamma>0$ and therefore has the positive real root $\lambda=(-\gamma)^{1/\alpha}$. Hence, the equilibrium is unstable.
		Now for $\tau_2 > 0$, define
		\begin{equation}
			F(\lambda) := \lambda^\alpha + \gamma - k e^{-\lambda \tau_1} + k e^{-(\lambda + \gamma)\tau_2 - \lambda \tau_1},
			\text{the characteristic function.}
			\label{CF}
		\end{equation}
		
		At $\lambda=0$,
		\[
		F(0) = \gamma - k + k e^{-\gamma \tau_2}.
		\]
		Since $\gamma<0$ and $\tau_2>0$, one has $e^{-\gamma\tau_2}>1$. Multiplication by $k<0$ reverses the inequality and gives $k e^{-\gamma\tau_2}<k$. Hence,
		\[
		F(0)=\gamma+\bigl(k e^{-\gamma\tau_2}-k\bigr)<\gamma<0.
		\]
		On the other hand, $\lim_{\lambda\to\infty}F(\lambda)=\infty$. The intermediate value theorem therefore guarantees a number $\lambda_*>0$ such that $F(\lambda_*)=0$. This positive characteristic root proves instability for every $\tau_1,\tau_2\geq0$.\\
		
		\textbf{Case (ii)} For $\alpha=1,$ the characteristic function (\ref{CF}) becomes 
		\[
		F(\lambda) := \lambda + \gamma - k e^{-\lambda \tau_1} + k e^{-(\lambda + \gamma)\tau_2 - \lambda \tau_1}.
		\]
		In the fourth quadrant, $\gamma<0$, so $-\gamma>0$. Direct substitution gives $F(-\gamma)=0$. Thus, $-\gamma$ is a positive real characteristic root, and the equilibrium is unstable.
	\end{proof}
	In the following theorem, we present a delay-dependent result for this case.
	\begin{theorem}
		\label{Case2Thrm2}
		In the fourth quadrant of the \(k\!-\!\gamma\) plane (so \(k>0\) and \(\gamma<0\)), fix \(0<\alpha<1\) and \(\tau_1,\tau_2>0\). Define the sufficient instability threshold
		\begin{equation}\label{eqkstar}
			k_* \;=
			\frac{(\tau_1+\tau_2)^{1+\frac{\tau_1}{\tau_2}}}{\tau_2 \; e^{\gamma\tau_1} \;\tau_1^{\tau_1 / \tau_2}}
			\left[\left(\frac{\log\left(1+\frac{\tau_2}{\tau_1}\right)}{\tau_2} - \gamma\right)^{\alpha} + \gamma \right],
		\end{equation}
		such that for every \(k>k_*\) the characteristic equation
		\[
		F(\lambda)=\lambda^\alpha+\gamma - k e^{-\lambda \tau_1} + k e^{-(\lambda+\gamma)\tau_2-\lambda \tau_1}
		\]
		has a root with positive real part; hence, the equilibrium is unstable for all such \(k\).
	\end{theorem}
	\begin{proof}
		Write
		\[
		F(\lambda)=L(\lambda)-R(\lambda),
		\]
		where
		\[
		L(\lambda)=\lambda^\alpha+\gamma
		\quad\text{and}\quad
		R(\lambda)=k e^{-\lambda\tau_1}\bigl(1-e^{-(\lambda+\gamma)\tau_2}\bigr).
		\]
		Since $L(0)=\gamma<0$ and $L'(\lambda)=\alpha\lambda^{\alpha-1}>0$ for $\lambda>0$, the function $L$ is strictly increasing on $(0,\infty)$. Moreover, $e^{-\gamma\tau_2}>1$ because $\gamma<0$, and hence
		\[
		R(0)=k\bigl(1-e^{-\gamma\tau_2}\bigr)<0.
		\]
		Also, $R(\lambda)\to0$ as $\lambda\to\infty$. To locate the global maximum of $R$, differentiate:
		\[
		R'(\lambda)=k e^{-\lambda \tau_1}H(\lambda),
		\qquad
		H(\lambda)=(\tau_1+\tau_2)e^{-(\lambda+\gamma)\tau_2}-\tau_1.
		\]
		The function $H$ is strictly decreasing because
		\[
		H'(\lambda)
		=-\tau_2(\tau_1+\tau_2)e^{-(\lambda+\gamma)\tau_2}<0.
		\]
		Moreover,
		\[
		H(0)=(\tau_1+\tau_2)e^{-\gamma\tau_2}-\tau_1
		>\tau_2>0,
		\]
		whereas $\lim_{\lambda\to\infty}H(\lambda)=-\tau_1<0$. Hence, $H$ has exactly one zero, namely
		\[
		\lambda_*=\frac{1}{\tau_2}\log\left(1+\frac{\tau_2}{\tau_1}\right)-\gamma>0.
		\]
		Since $k e^{-\lambda\tau_1}>0$, the sign of $R'(\lambda)$ is the sign of $H(\lambda)$. Thus, $R$ is strictly increasing on $(0,\lambda_*)$ and strictly decreasing on $(\lambda_*,\infty)$. Therefore, $\lambda_*$ is the unique global maximizer of $R$.

		At $\lambda=\lambda_*$,
		\[
		e^{-(\lambda_*+\gamma)\tau_2}
		=\frac{\tau_1}{\tau_1+\tau_2}
		\]
		and
		\[
		e^{-\lambda_*\tau_1}
		=e^{\gamma\tau_1}
		\left(\frac{\tau_1}{\tau_1+\tau_2}\right)^{\tau_1/\tau_2}.
		\]
		Consequently,
		\[
		R(\lambda_*)=kA,
		\qquad
		A=
		\frac{\tau_2 e^{\gamma\tau_1}\tau_1^{\tau_1/\tau_2}}
		{(\tau_1+\tau_2)^{1+\tau_1/\tau_2}}>0.
		\]
		Define $k_*$ by $k_*A=L(\lambda_*)$. Equivalently,
		\[
		\lambda_*^\alpha+\gamma
		=k_*e^{-\lambda_*\tau_1}
		\bigl(1-e^{-(\lambda_*+\gamma)\tau_2}\bigr).
		\]
		Substitution of the preceding expressions gives
		\[
		\left(\frac{1}{\tau_2}\log\left(1+\frac{\tau_2}{\tau_1}\right)-\gamma\right)^\alpha+\gamma
		=
		\frac{k_*\tau_2\tau_1^{\tau_1/\tau_2}e^{\gamma\tau_1}}
		{(\tau_1+\tau_2)^{1+\tau_1/\tau_2}}.
		\]
		Solving this equality for $k_*$ yields
		\[
		k_*=
		\frac{(\tau_1+\tau_2)^{1+\frac{\tau_1}{\tau_2}}}
		{\tau_2 e^{\gamma\tau_1}\tau_1^{\tau_1/\tau_2}}
		\left[
		\left(\frac{\log\left(1+\frac{\tau_2}{\tau_1}\right)}{\tau_2}-\gamma\right)^\alpha+\gamma
		\right].
		\]
		
		\begin{figure}[H]
			\centering
			\includegraphics[width=0.65\textwidth]{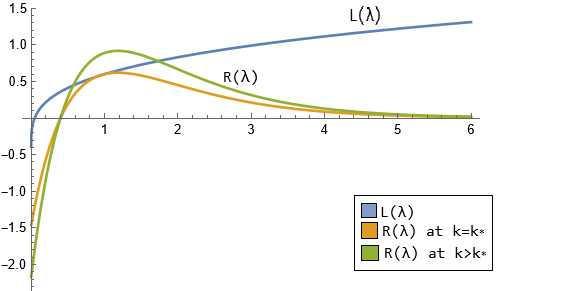}
			\caption{Comparison of $L(\lambda)=\lambda^\alpha+\gamma$ and $R(\lambda)$.
				At $k=k_*$ the curves intersect at $\lambda_*$, where $R$ attains its global maximum.
				For $k>k_*$, one has $R(\lambda_*)>L(\lambda_*)$, which yields a positive root of $F(\lambda)$ and hence instability.}
			\label{figtheorem_illustration}
		\end{figure}
		We now prove that $k>k_*$ is sufficient for instability. Since $A>0$,
		\[
		k>k_*
		\quad\Longrightarrow\quad
		R(\lambda_*)=kA>k_*A=L(\lambda_*).
		\]
		Therefore,
		\[
		F(\lambda_*)=L(\lambda_*)-R(\lambda_*)<0.
		\]
		On the other hand, the exponential terms in $F$ tend to zero as $\lambda\to\infty$, while $\lambda^\alpha\to\infty$. Hence,
		\[
		\lim_{\lambda\to\infty}F(\lambda)=\infty.
		\]
		There is therefore a number $\Lambda>\lambda_*$ for which $F(\Lambda)>0$. Since $F$ is continuous on $[\lambda_*,\Lambda]$, the intermediate value theorem gives a point $\lambda_0\in(\lambda_*,\Lambda)$ satisfying
		\[
		F(\lambda_0)=0.
		\]
		The number $\lambda_0$ is a positive real characteristic root. A positive characteristic root produces an exponentially growing component of the linearized solution and rules out asymptotic stability. Thus, the equilibrium is unstable whenever $k>k_*$. This argument proves sufficiency only; it does not claim that $k_*$ is a necessary instability threshold.
	\end{proof}

	{
	\subsection*{Hopf interpretation for two nonzero delays}
	When both delays are retained, a Hopf candidate must satisfy the real and imaginary equations
	\eqref{eqreal_case2}--\eqref{eqimag_case2} for some $v>0$. A Hopf bifurcation occurs only if the corresponding pair $\lambda=\pm iv$ is simple and crosses the imaginary axis with nonzero speed as the selected delay parameter varies. The proof of Theorem~\ref{thmcase2_unified}(a) shows that, in the regions
	\[
	\gamma>2k>0
	\qquad\text{or}\qquad
	\gamma>-2k>0,
	\]
	the real crossing equation cannot be satisfied. These delay-independent stability regions therefore contain no Hopf bifurcation for any $\tau_1,\tau_2\geq0$.

	In parameter regions not excluded by Theorem~\ref{thmcase2_unified}, Equations~\eqref{eqreal_case2} and~\eqref{eqimag_case2}, together with a nonzero transversality derivative, provide the direct extension of the Hopf condition in Theorem~\ref{thmprelim}. The instability mechanism in Theorem~\ref{Case2Thrm2} is different: it is proved by the existence of a positive real characteristic root when $k>k_*$. Thus, $k=k_*$ is a real-root instability threshold and should not be identified as a Hopf boundary. This distinction connects the stability classifications to the preliminary Hopf theorem and separates oscillatory loss of stability from instability caused by a positive real root.
	}

	\begin{figure}[H]
		\centering
		\includegraphics[scale=0.58]{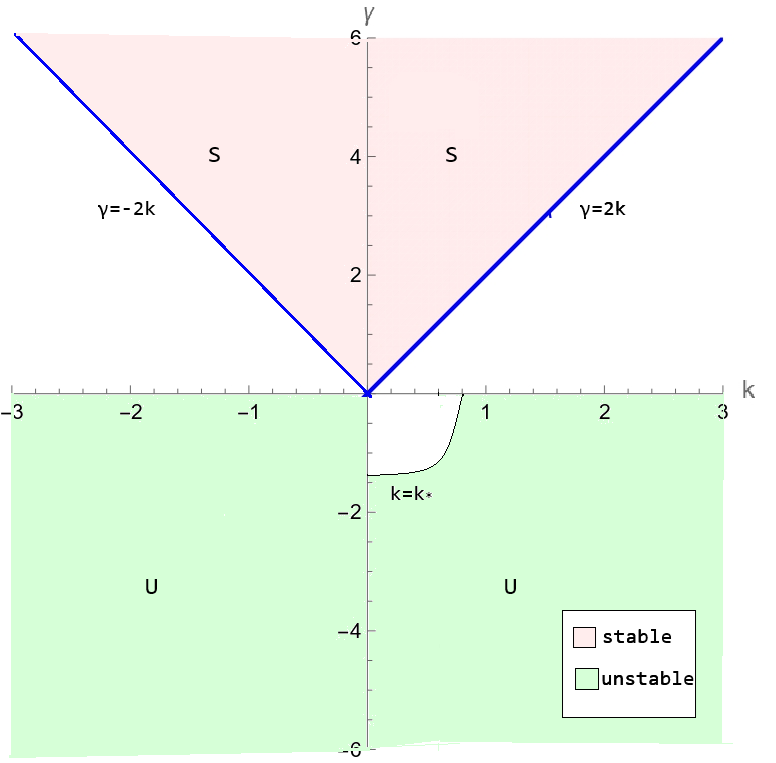}
		\caption{Stability regions for $\tau_1>0$}
		\label{ind2}
	\end{figure}
	\subsection{Illustrative Examples for Case 2}
	Consider the FDDE $D^{\alpha} x(t) = -\gamma x(t) +kx(t-\tau_1)-k e^{-\gamma\tau_2} x(t - \tau_1-\tau_2)$.\\
    
	\textbf{Example 5.1.1}
	
\begin{figure}[h!]
    \centering
    \begin{subfigure}{0.45\textwidth}
        \centering
        \includegraphics[width=\linewidth]{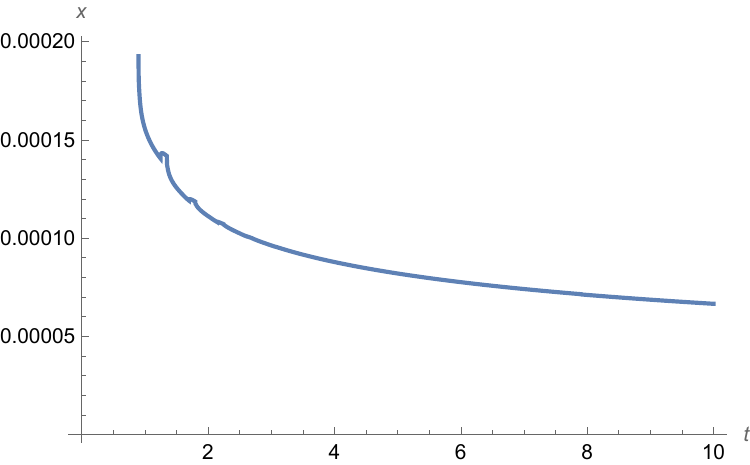}
        \caption{$\alpha=0.33,\; k=2.6,\; \gamma=6.4,\; \tau=0.35,\;\tau_2=0.45$}
        \label{5a}
    \end{subfigure}
    \hfill
    \begin{subfigure}{0.45\textwidth}
        \centering
        \includegraphics[width=\linewidth]{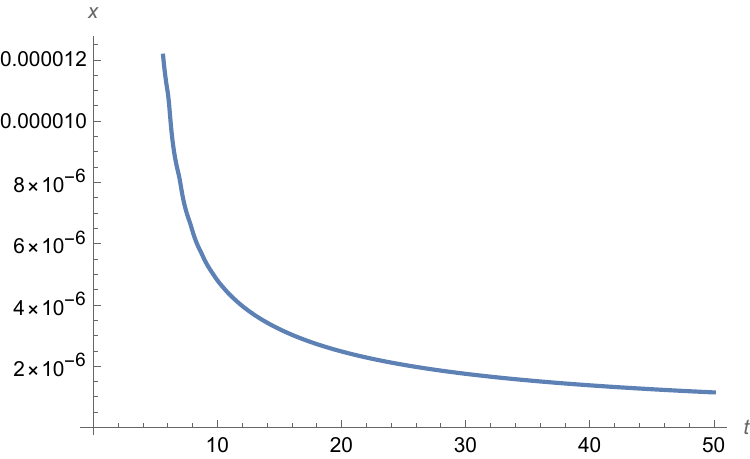}
        \caption{$\alpha=0.8,\; k=6.4,\; \gamma=15,\; \tau=0.85,\;\tau_2=1.45$}
        \label{5b}
    \end{subfigure}
    \caption{$\gamma>2k>0$}
\end{figure}

	Take $\alpha=0.33,\; k=2.6,\; \gamma=6.4,\; \tau_1=0.35,\;\tau_2=0.45$.  
	Since $\gamma>2k>0$, Theorem~\ref{thmcase2_unified}(a)(i) guarantees delay-independent stability. Simulations confirm convergence (see Figure (\ref{5a})).\\

    	\textbf{Example 5.1.2}
	
	Take $\alpha=0.8,\; k=6.4,\; \gamma=15,\; \tau_1=0.85,\;\tau_2=1.45$.  
Again, $\gamma>2k>0$. Theorem~\ref{thmcase2_unified}(a)(i) therefore guarantees stability for these parameter values; see Figure~\ref{5b}.\\
	
	\textbf{Example 5.1.3}
    
\begin{figure}[h!]
    \centering
    \begin{subfigure}{0.45\textwidth}
        \centering
        \includegraphics[width=\linewidth]{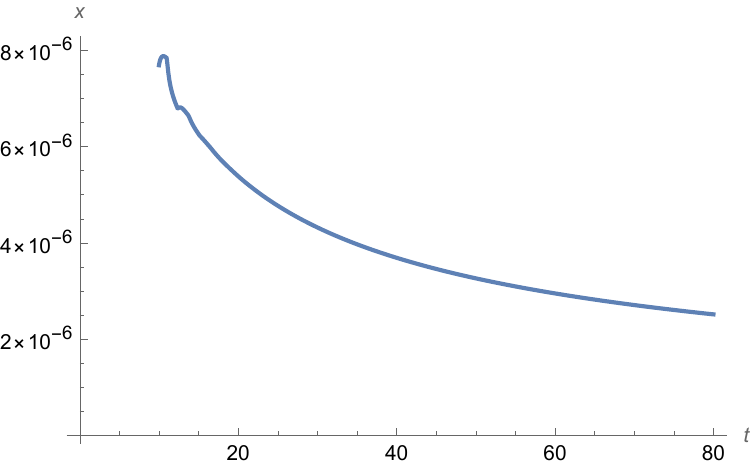}
        \caption{$\alpha=0.56,\; k=-4.9,\; \gamma=12,\; \tau_1=1.37,\; \tau_2=1.49$}
        \label{5c}
    \end{subfigure}
    \hfill
    \begin{subfigure}{0.45\textwidth}
        \centering
        \includegraphics[width=\linewidth]{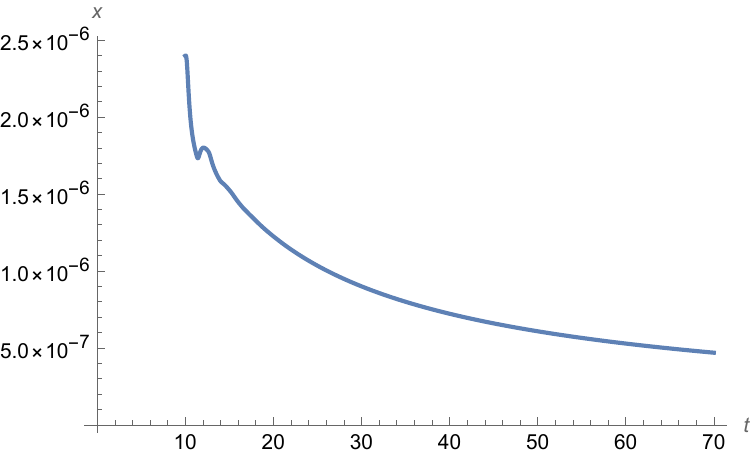}
        \caption{$\alpha=0.78,\; k=-5.2,\; \gamma=13.4,\; \tau_1=1.25,\; \tau_2=0.64$}
        \label{5d}
    \end{subfigure}
    \caption{$\gamma>-2k>0$}
\end{figure}

	Now, consider $\alpha=0.56,\; k=-4.9,\; \gamma=12,\; \tau_1=1.37,\; \tau_2=1.49$.  
	Here $\gamma>-2k>0$. By Theorem~\ref{thmcase2_unified}(a)(ii), the system is asymptotically stable for all $\tau$. Numerical solutions confirm this (cf. Figure (\ref{5c})).\\

    	\textbf{Example 5.1.4}
	
	Take  $\alpha=0.78,\; k=-5.2,\; \gamma=13.4,\; \tau_1=1.25,\; \tau_2=0.64$.  
Again, $\gamma>-2k>0$, so Theorem~\ref{thmcase2_unified}(a)(ii) guarantees stability; see Figure~\ref{5d}.\\
	
	\textbf{Example 5.1.5}
    
\begin{figure}[h!]
    \centering
    \begin{subfigure}{0.45\textwidth}
        \centering
        \includegraphics[width=\linewidth]{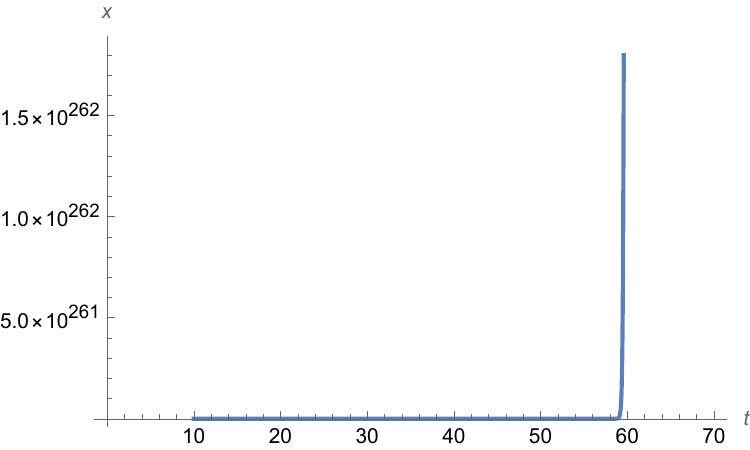}
        \caption{ $\alpha=0.65,\; k=-2.35,\; \gamma=-4.56,\; \tau_1=0.82,\; \tau_2=0.41$}
        \label{5e}
    \end{subfigure}
    \hfill
    \begin{subfigure}{0.45\textwidth}
        \centering
        \includegraphics[width=\linewidth]{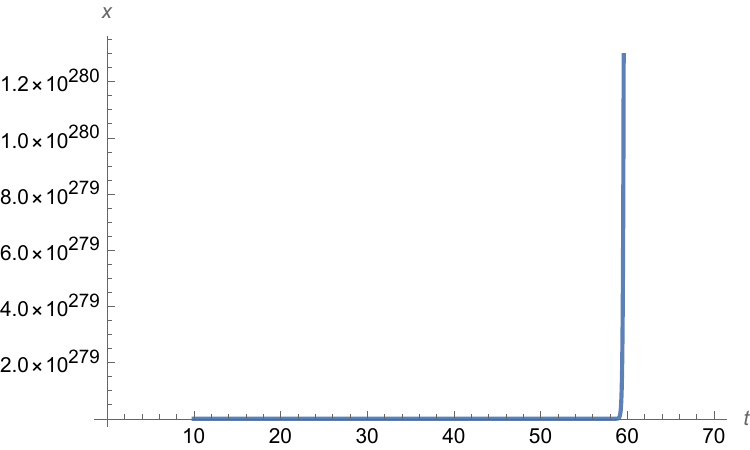}
        \caption{$\alpha=0.89,\; k=-3.31,\; \gamma=-8.41,\; \tau_1=0.42,\; \tau_2=0.62$}
        \label{5f}
    \end{subfigure}
    \caption{$k<0,\;\gamma<0$}
\end{figure}

	For $\alpha=0.65,\; k=-2.35,\; \gamma=-4.56,\; \tau_1=0.82,\; \tau_2=0.41$.  
	Since $k<0,\;\gamma<0$, Theorem~\ref{thmcase2_unified}(b) predicts delay-independent instability. The simulations show divergence  (see Figure (\ref{5e})).\\

    \textbf{Example 5.1.6}
    
Now, take $\alpha=0.89,\; k=-3.31,\; \gamma=-8.41,\; \tau_1=0.42,\; \tau_2=0.62$.  
	Here $(k,\gamma)$ again lies in the third quadrant. Theorem~\ref{thmcase2_unified}(b) therefore places the equilibrium in the unstable region; see Figure~\ref{5f}.\\
	\textbf{Note:} By testing a range of $\tau$ values, we verify that stability is independent of delay; however, for the sake of conciseness, we only show a few representative examples here.\\
	
	\textbf{Example 5.1.7}
    
	For $\alpha=0.43,\; k=14,\; \gamma=-2.32,\; \tau_1=1.28,\; \tau_2=2.35$.  
	Since $k>0$ and $\gamma<0$, Theorem~\ref{Case2Thrm2} provides a threshold $k_*$ above which the equilibrium is unstable. For the present parameters, $k_*=11.3536$ and $k=14>k_*$. The equilibrium is therefore unstable.
	\begin{figure}[H]
		\centering
		\includegraphics[scale=0.5]{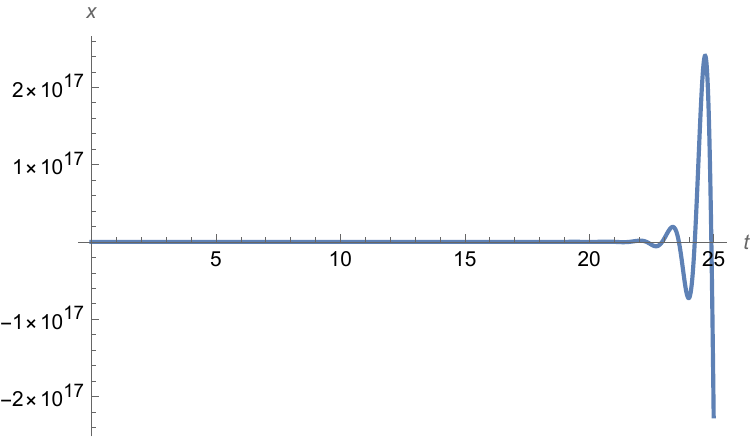}
		\caption{$\alpha=0.43,\; k=14,\; \gamma=-2.32,\; \tau_1=1.28,\; \tau_2=2.35$}
        \label{5g}
	\end{figure}
	
	\subsection{Illustrative Examples for a Nonlinear Function $g(x)$}
	
	Let us consider the nonlinear fractional delay differential equation (\ref{eqmain}) with  
	\[
	g(x) = k \sin x,
	\]
	which leads to
	\begin{equation}
		\label{nleg}
		D^\alpha x(t) = -\gamma x(t) + k \sin\!\bigl(x(t-\tau_1)\bigr) 
		- e^{-\gamma \tau_2} \, k \sin\!\bigl(x(t-\tau_1-\tau_2)\bigr), 
		\quad 0 < \alpha \leq 1.
	\end{equation}
	
	To analyze the stability of the trivial equilibrium $x_* = 0$, we linearize the system around this point.  
	Since
	\[
	g'(x) = k \cos x \quad \text{and therefore} \quad g'(0) = k,
	\]
	the linearized equation near $x_* = 0$ takes the form
	\[
	D^\alpha x(t) = -\gamma x(t) + k x(t-\tau_1) - k e^{-\gamma \tau_2} x(t-\tau_1-\tau_2).
	\]
	
\begin{figure}[H]
		\centering
		\includegraphics[scale=0.5]{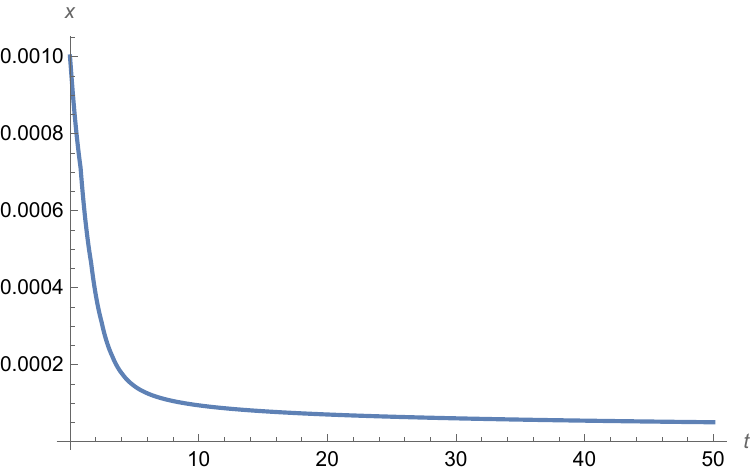}
		\caption{$\alpha = 0.35,\; k = 2.4,\; \gamma = 6,\; \tau_1 = 0.48,\; \tau_2=1.29$.}
		\label{fig5h}
	\end{figure}

	This is exactly of the same type as the linearized equation considered in Section~\ref{secmain}. Hence, the stability of the trivial solution is again determined by the conditions given in Theorem~\ref{thmcase2_unified}.  
	We now illustrate this with a few examples.\\
	\textbf{Example 5.2.1.}  
	Let $\alpha = 0.35,\; k = 2.4,\; \gamma = 6,\; \tau_1 = 0.82,\; \tau_2=1.49$.  
	Then, equation~\eqref{nleg} becomes
	\[
	D^{0.35} x(t) = -6x(t) + 2.4 \sin(t-0.48) - 2.4 e^{-8.94} \sin(t-2.32).
	\]
	The solutions are stable (cf. Figure \ref{fig5h}), in agreement with Theorem~\ref{thmcase2_unified}, since here $\gamma > 2k > 0$.  
	
\begin{figure}[H]
		
		\centering
		\includegraphics[scale=0.5]{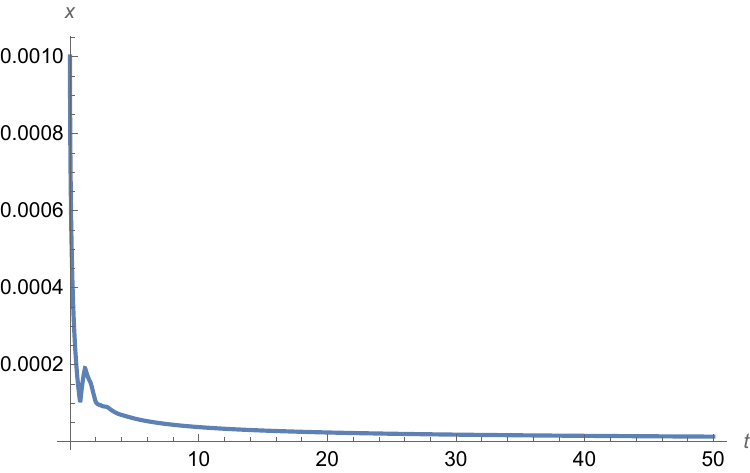}
		\caption{$\alpha = 0.65,\; k = -0.8,\; \gamma = 2,\; \tau_1 = 0.79,\; \tau_2=0.38$.}
		\label{fig5i}
	\end{figure}

	\textbf{Example 5.2.2.}  
	Now take $\alpha = 0.65,\; k = -0.8,\; \gamma = 2,\; \tau_1 = 0.79,\; \tau_2=0.38$.  
	In this case,
	\[
	D^{0.65} x(t) = -2x(t) - 0.8 \sin(t-0.79) + 0.8 e^{-0.76} \sin(t-1.17).
	\]
	The solutions converge to the equilibrium (Figure \ref{fig5i}), again supporting Theorem~\ref{thmcase2_unified}, since $\gamma > -2k > 0$.  
	
\begin{figure}[H]
		
		\centering
		\includegraphics[scale=0.5]{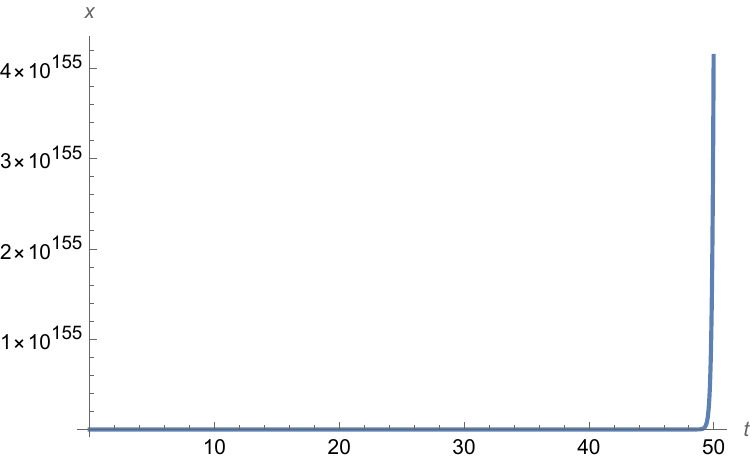}
		\caption{$\alpha = 0.88,\; k = -3.5,\; \gamma = -5.7,\; \tau_1 = 1.62,\; \tau_2=0.85$.}
		\label{fig5j}
	\end{figure}

	\textbf{Example 5.2.3.}  
	Finally, consider $\alpha = 0.88,\; k = -3.5,\; \gamma = -5.7,\;  \tau_1 = 1.62,\; \tau_2=0.85$.  
	Equation~\eqref{nleg} becomes
	\[
	D^{0.88} x(t) = 5.7x(t) - 3.5 \sin(t-1.62) + 3.5 e^{4.845} \sin(t-2.47).
	\]
	Here, the solutions are unstable (Figure \ref{fig5j}), consistent with Theorem~\ref{thmcase2_unified}, since $(k,\gamma)$ lies in the third quadrant.  
	
	\section{Conclusion}
\label{secconclusion}
{
We studied a nonlinear Caputo fractional delay differential equation in which two discrete delays interact with the delay-dependent gain $e^{-\gamma\tau_2}$. Before performing the stability analysis, we established local existence and uniqueness under a local Lipschitz assumption on $g$. When $g$ is globally Lipschitz, the unique solution exists for all $t\geq0$.

For $\tau_1=0$, the linearized equation becomes a single-delay problem with the nonconstant coefficient $-k e^{-\gamma\tau}$. Tracking this coefficient through the known stability regions yields delay-independent stable and unstable parameter sets and identifies the implicit critical-delay boundary $\tau=\tau_1^*(\tau)$. Subject to simplicity and transversality of the imaginary roots, this boundary represents a Hopf bifurcation, whereas $\tau_{2*}$ only marks a transition between coefficient regions. For $\tau_1>0$ and $\tau_2\geq0$, careful estimates of the characteristic equation exclude imaginary-axis roots in the stated delay-independent stability regions. We also proved delay-independent instability in the third quadrant and, for $0<\alpha<1$ in the fourth quadrant, obtained the explicit sufficient threshold $k_*$ that guarantees a positive real characteristic root.

The analysis demonstrates that the fractional order affects both the spectral geometry and the critical conditions through the factors $\cos(\alpha\pi/2)$, $\sin(\alpha\pi/2)$, and $\lambda^\alpha$. Consequently, the fractional problem cannot be reduced to its integer-order counterpart by a direct parameter substitution. The numerical examples and stability diagrams are consistent with the analytical classifications.

The results are primarily local because they are based on linearization, and several conditions for the two-delay case are sufficient rather than necessary. A complete characterization of the remaining parameter regions, rigorous nonlinear Hopf analysis, normal-form calculations, and parameter estimation for biological data are natural directions for future work.
}

\section*{Acknowledgment}	
	Pragati Dutta thanks the University of Hyderabad for the non-net fellowship (UH/DSW/F\&/Non-NET/2023/165).
	 The authors are thankful to the anonymous referee for their insightful comments, which led to an improved manuscript.
	 
\section*{Declarations}
The authors declare that they have no known competing financial interests or personal relationships that could have appeared to influence the work reported in this paper.\\
Generative AI tools were used only to improve the language and 
clarity of this manuscript. The authors carefully reviewed and 
edited the output and take full responsibility for the final content.
		\bibliographystyle{unsrt}	
	\bibliography{reff}

@book{podlubny1998fractional,
  title={Fractional differential equations: An introduction to fractional derivatives, fractional differential equations, to methods of their solution and some of their applications},
  author={Podlubny, Igor},
  volume={198},
  year={1998},
  publisher={elsevier}
}

@article{boullu2020stability,
	title={Stability analysis of an equation with two delays and application to the production of platelets},
	author={Boullu, Lo{\"\i}s and Pujo-Menjouet, Laurent and B{\'e}lair, Jacques},
	journal={Discrete and Continuous Dynamical Systems-Series S},
	volume={13},
	number={11},
	pages={3005--3027},
	year={2020}
}

@article{belair1994stability,
	title={Stability and bifurcations of equilibria in a multiple-delayed differential equation},
	author={B{\'e}lair, Jacques and Campbell, Sue Ann},
	journal={SIAM Journal on Applied Mathematics},
	volume={54},
	number={5},
	pages={1402--1424},
	year={1994},
	publisher={SIAM}
}

@article{piotrowska2008hopf,
	title={Hopf bifurcation in a solid avascular tumour growth model with two discrete delays},
	author={Piotrowska, Monika Joanna},
	journal={Mathematical and Computer Modelling},
	volume={47},
	number={5-6},
	pages={597--603},
	year={2008},
	publisher={Elsevier}
}

@book{kilbas2006theory,
  title={Theory and applications of fractional differential equations},
  author={Kilbas, Anatoli{\u\i} Aleksandrovich and Srivastava, Hari M and Trujillo, Juan J},
  volume={204},
  year={2006},
  publisher={elsevier}
}

@book{diethelm2010analysis,
  title={The Analysis of Fractional Differential Equations: An Application-Oriented Exposition Using Differential Operators of Caputo Type},
  author={Diethelm, Kai},
  series={Lecture Notes in Mathematics},
  volume={2004},
  publisher={Springer},
  address={Berlin, Heidelberg},
  year={2010},
  doi={10.1007/978-3-642-14574-2}
}

@article{ye2007generalized,
  title={A generalized {Gronwall} inequality and its application to a fractional differential equation},
  author={Ye, Haiping and Gao, Jianming and Ding, Yongsheng},
  journal={Journal of Mathematical Analysis and Applications},
  volume={328},
  number={2},
  pages={1075--1081},
  year={2007},
  doi={10.1016/j.jmaa.2006.05.061}
}

@article{ye2011henry,
  title={Henry--Gronwall type retarded integral inequalities and their applications to fractional differential equations with delay},
  author={Ye, Haiping and Gao, Jianming},
  journal={Applied Mathematics and Computation},
  volume={218},
  number={8},
  pages={4152--4160},
  year={2011},
  doi={10.1016/j.amc.2011.09.046}
}

@book{hale2013introduction,
  title={Introduction to functional differential equations},
  author={Hale, Jack K and Lunel, Sjoerd M Verduyn},
  volume={99},
  year={2013},
  publisher={Springer Science \& Business Media}
}

@article{bhalekar2016stability,
  title={Stability and bifurcation analysis of a generalized scalar delay differential equation},
  author={Bhalekar, Sachin},
  journal={Chaos: An Interdisciplinary Journal of Nonlinear Science},
  volume={26},
  number={8},
  year={2016},
  publisher={AIP Publishing}
}

@article{bhalekar2025analysis,
  title={Analysis of a class of two-delay fractional differential equation},
  author={Bhalekar, Sachin and Dutta, Pragati},
  journal={Chaos: An Interdisciplinary Journal of Nonlinear Science},
  volume={35},
  number={1},
  year={2025},
  publisher={AIP Publishing}
}

@article{gupta2024fractional,
  title={Fractional order sunflower equation: stability, bifurcation and chaos},
  author={Gupta, Deepa and Bhalekar, Sachin},
  journal={The European Physical Journal Special Topics},
  pages={1--11},
  year={2024},
  publisher={Springer}
}

@book{mainardi2022fractional,
  title={Fractional calculus and waves in linear viscoelasticity: an introduction to mathematical models},
  author={Mainardi, Francesco},
  year={2022},
  publisher={World Scientific}
}

@article{li1999stability,
  title={Stability and bifurcation in delay--differential equations with two delays},
  author={Li, Xiangao and Ruan, Shigui and Wei, Junjie},
  journal={Journal of Mathematical Analysis and Applications},
  volume={236},
  number={2},
  pages={254--280},
  year={1999},
  publisher={Elsevier}
}

@article{hale1993global,
  title={Global geometry of the stable regions for two delay differential equations},
  author={Hale, Jack K and Huang, WZ},
  journal={Journal of Mathematical analysis and applications},
  volume={178},
  number={2},
  pages={344--362},
  year={1993},
  publisher={Academic Press}
}

@article{daftardar2012dynamics,
  title={Dynamics of fractional-ordered {Chen} system with delay},
  author={Daftardar-Gejji, Varsha and Bhalekar, Sachin and Gade, Prashant},
  journal={Pramana},
  volume={79},
  pages={61--69},
  year={2012},
  publisher={Springer}
}

@article{bhalekar2011fractional,
  title={Fractional Bloch equation with delay},
  author={Bhalekar, Sachin and Daftardar-Gejji, Varsha and Baleanu, Dumitru and Magin, Richard},
  journal={Computers \& Mathematics with Applications},
  volume={61},
  number={5},
  pages={1355--1365},
  year={2011},
  publisher={Elsevier}
}

@article{bhalekar2012dynamical,
  title={Dynamical analysis of fractional order U{\c{c}}ar prototype delayed system},
  author={Bhalekar, Sachin},
  journal={Signal, Image and Video Processing},
  volume={6},
  number={3},
  pages={513--519},
  year={2012},
  publisher={Springer}
}

@article{bhalekar2012generalized,
  title={Generalized fractional order Bloch equation with extended delay},
  author={Bhalekar, Sachin and Daftardar-Gejji, Varsha and Baleanu, Dumitru and Magin, Richard},
  journal={International Journal of Bifurcation and Chaos},
  volume={22},
  number={04},
  pages={1250071},
  year={2012},
  publisher={World Scientific}
}

@book{Podlubny1999,
  author    = {Igor Podlubny},
  title     = {Fractional Differential Equations},
  publisher = {Academic Press},
  year      = {1999},
  address   = {San Diego},
}

@book{smith2010introduction,
  title={An introduction to delay differential equations with applications to the life sciences},
  author={Smith, Hal L},
  year={2010},
  publisher={springer}
}

@book{lakshmanan2011dynamics,
  title={Dynamics of nonlinear time-delay systems},
  author={Lakshmanan, Muthusamy and Senthilkumar, Dharmapuri Vijayan},
  year={2011},
  publisher={Springer Science \& Business Media}
}

@article{epstein1991differential,
  title={Differential delay equations in chemical kinetics. Nonlinear models: The cross-shaped phase diagram and the Oregonator},
  author={Epstein, Irving R and Luo, Yin},
  journal={The Journal of chemical physics},
  volume={95},
  number={1},
  pages={244--254},
  year={1991},
  publisher={American Institute of Physics}
}

@book{rihan2021delay,
  title={Delay Differential Equations and Applications to Biology},
  author={Rihan, Fathalla A},
  year={2021},
  publisher={Springer Nature}
}

@article{bani2017analysis,
  title={Analysis and applications of delay differential equations in biology and medicine},
  author={Bani-Yaghoub, Majid},
  journal={arXiv preprint arXiv:1701.04173},
  year={2017}
}

@article{roussel1996use,
  title={The use of delay differential equations in chemical kinetics},
  author={Roussel, Marc R},
  journal={The journal of physical chemistry},
  volume={100},
  number={20},
  pages={8323--8330},
  year={1996},
  publisher={ACS Publications}
}

@article{flore2019feynman,
  title={A {Feynman-Kac} type formula for a fixed delay {CIR} model},
  author={Flore, Federico and Nappo, Giovanna},
  journal={Stochastic Analysis and Applications},
  volume={37},
  number={4},
  pages={550--573},
  year={2019},
  publisher={Taylor \& Francis}
}

@article{agrawal2020jump,
  title={Jump models with delay—option pricing and logarithmic {Euler--Maruyama} scheme},
  author={Agrawal, Nishant and Hu, Yaozhong},
  journal={Mathematics},
  volume={8},
  number={11},
  pages={1932},
  year={2020},
  publisher={MDPI}
}

@article{bhalekar2011predictor,
  title={A predictor-corrector scheme for solving nonlinear delay differential equations of fractional order},
  author={Bhalekar, Sachin and Daftardar-Gejji, Varsha},
  journal={J. Fract. Calc. Appl},
  volume={1},
  number={5},
  pages={1--9},
  year={2011}
}

@article{daftardar2014new,
  title={A new predictor--corrector method for fractional differential equations},
  author={Daftardar-Gejji, Varsha and Sukale, Yogita and Bhalekar, Sachin},
  journal={Applied Mathematics and Computation},
  volume={244},
  pages={158--182},
  year={2014},
  publisher={Elsevier}
}

@article{braddock1983two,
  title={On a two lag differential delay equation},
  author={Braddock, RD and Van den Driessche, P},
  journal={The ANZIAM Journal},
  volume={24},
  number={3},
  pages={292--317},
  year={1983},
  publisher={Cambridge University Press}
}

@article{bhalekar2019analysing,
  title={Analysing the stability of a delay differential equation involving two delays},
  author={Bhalekar, Sachin},
  journal={Pramana},
  volume={93},
  number={2},
  pages={24},
  year={2019},
  publisher={Springer}
}

@article{sprott2007simple,
  title={A simple chaotic delay differential equation},
  author={Sprott, JC},
  journal={Physics Letters A},
  volume={366},
  number={4-5},
  pages={397--402},
  year={2007},
  publisher={Elsevier}
}

@article{senthilkumar2005bifurcations,
  title={Bifurcations and chaos in time delayed piecewise linear dynamical systems},
  author={Senthilkumar, DV and Lakshmanan, M},
  journal={International Journal of Bifurcation and Chaos},
  volume={15},
  number={09},
  pages={2895--2912},
  year={2005},
  publisher={World Scientific}
}

@article{ruiz2013chaos,
  title={Chaos in delay differential equations with applications in population dynamics},
  author={Ruiz-Herrera, Alfonso},
  journal={Discrete Contin. Dyn. Syst},
  volume={33},
  number={4},
  pages={1633--1644},
  year={2013}
}

@article{glass2021nonlinear,
  title={Nonlinear delay differential equations and their application to modeling biological network motifs},
  author={Glass, David S and Jin, Xiaofan and Riedel-Kruse, Ingmar H},
  journal={Nature communications},
  volume={12},
  number={1},
  pages={1788},
  year={2021},
  publisher={Nature Publishing Group UK London}
}

@article{otto2019nonlinear,
  title={Nonlinear dynamics of delay systems: an overview},
  author={Otto, A and Just, W and Radons, G},
  journal={Philosophical Transactions of the Royal Society A},
  volume={377},
  number={2153},
  pages={20180389},
  year={2019},
  publisher={The Royal Society Publishing}
}

@article{bazighifan2021oscillation,
  title={On the oscillation of nonlinear delay differential equations and their applications},
  author={Bazighifan, Omar and Askar, Sameh},
  journal={Open Physics},
  volume={19},
  number={1},
  pages={788--796},
  year={2021},
  publisher={De Gruyter}
}

@incollection{gumussoy2014computer,
  title={Computer Aided Control System Design for Time Delay Systems Using MATLAB{\textregistered}},
  author={Gumussoy, Suat and Gahinet, Pascal},
  booktitle={Delay Systems: From Theory to Numerics and Applications},
  pages={257--270},
  year={2014},
  publisher={Springer}
}

@article{nagarajan2025stock,
  title={Stock Return Model Using Stochastic Delay Differential Equation in Finance.},
  author={Nagarajan, Racshitha and Rajendran, Manimaran and Chandrasekaran, Vijayan},
  journal={Mathematical Modelling of Engineering Problems},
  volume={12},
  number={1},
  year={2025}
}

@article{belair1995age,
  title={Age-structured and two-delay models for erythropoiesis},
  author={B{\'e}lair, Jacques and Mackey, Michael C and Mahaffy, Joseph M},
  journal={Mathematical biosciences},
  volume={128},
  number={1-2},
  pages={317--346},
  year={1995},
  publisher={Elsevier}
}

@book{gu2003stability,
  title={Stability of time-delay systems},
  author={Gu, Keqin and Chen, Jie and Kharitonov, Vladimir L},
  year={2003},
  publisher={Springer Science \& Business Media}
}

@book{niculescu2002delay,
  title={Delay effects on stability: a robust control approach},
  author={Niculescu, Silviu-Iulian},
  year={2002},
  publisher={Springer}
}

@article{bhalekar2019can,
  title={Can we split fractional derivative while analyzing fractional differential equations?},
  author={Bhalekar, Sachin and Patil, Madhuri},
  journal={Communications in Nonlinear Science and Numerical Simulation},
  volume={76},
  pages={12--24},
  year={2019},
  publisher={Elsevier}
}

@article{li2007remarks,
  title={Remarks on fractional derivatives},
  author={Li, Changpin and Deng, Weihua},
  journal={Applied mathematics and computation},
  volume={187},
  number={2},
  pages={777--784},
  year={2007},
  publisher={Elsevier}
}

@article{bagley1983theoretical,
  title={A theoretical basis for the application of fractional calculus to viscoelasticity},
  author={Bagley, Ronald L and Torvik, Peter J},
  journal={Journal of rheology},
  volume={27},
  number={3},
  pages={201--210},
  year={1983},
  publisher={The Society of Rheology}
}

@article{dong2023application,
  title={Application of a time-delay SIR model with vaccination in COVID-19 prediction and its optimal control strategy},
  author={Dong, Suyalatu and Xu, Linlin and A, Yana and Lan, Zhong-Zhou and Xiao, Ding and Gao, Bo},
  journal={Nonlinear Dynamics},
  volume={111},
  number={11},
  pages={10677--10692},
  year={2023},
  publisher={Springer}
}

@article{bhalekar2022stability,
  title={Stability and bifurcation analysis of a fractional order delay differential equation involving cubic nonlinearity},
  author={Bhalekar, Sachin and Gupta, Deepa},
  journal={Chaos, Solitons \& Fractals},
  volume={162},
  pages={112483},
  year={2022},
  publisher={Elsevier}
}

@book{insperger2011semi,
  title={Semi-discretization for time-delay systems: stability and engineering applications},
  author={Insperger, Tam{\'a}s and St{\'e}p{\'a}n, G{\'a}bor},
  volume={178},
  year={2011},
  publisher={Springer Science \& Business Media}
}

@article{michiels2010control,
  title={Control design for time-delay systems based on quasi-direct pole placement},
  author={Michiels, Wim and Vyhl{\'\i}dal, Tom{\'a}{\v{s}} and Z{\'\i}tek, Pavel},
  journal={Journal of Process Control},
  volume={20},
  number={3},
  pages={337--343},
  year={2010},
  publisher={Elsevier}
}

@book{leibniz1860leibnizens,
  title={Leibnizens mathematische schriften},
  author={Leibniz, Gottfried Wilhelm},
  volume={2},
  year={1860},
  publisher={Schmidt}
}

@article{deng2007stability,
  title={Stability analysis of linear fractional differential system with multiple time delays},
  author={Deng, Weihua and Li, Changpin and Lü, Jinhu},
  journal={Nonlinear Dynamics},
  volume={48},
  number={4},
  pages={409--416},
  year={2007},
  doi={10.1007/s11071-006-9094-0}
}

@article{abdalla2026existence,
  title={Existence and stability analysis of impulsive {Caputo} fractional delay systems},
  author={Abdalla, Manal Elzain Mohamed and Liu, Zhenhai and Hammad, Hasanen A.},
  journal={Boundary Value Problems},
  volume={2026},
  number={1},
  pages={32},
  year={2026},
  doi={10.1186/s13661-025-02214-4}
}

@article{elsaka2026dynamic,
  title={Dynamic analysis of the fractional distributed delay models},
  author={El-Saka, H. A. A. and El-Sherbeny, D. El. A. and El-Sayed, A. M. A.},
  journal={Scientific Reports},
  volume={16},
  pages={16252},
  year={2026},
  doi={10.1038/s41598-026-52327-8}
}

@article{zhu2026hopf,
  title={Hopf bifurcation and stability analysis of a fractional-order {Lotka--Volterra} predator--prey model with two delays},
  author={Zhu, Mengfan and Cheng, Zunshui and Xin, Youming and Shang, Yun and Lin, Xue and Cao, Jinde},
  journal={Intelligence \& Control},
  volume={2},
  number={2},
  pages={2},
  year={2026},
  doi={10.53941/ic.2026.100005}
}
	
\end{document}